\documentclass{amsart} \usepackage{amscd,amsmath,amssymb}
 \usepackage{epsfig} \usepackage{xypic} \usepackage{graphicx}
 \newcommand{\Alt}{\mathfrak{A}} \newcommand{\Vier}{\mathfrak{V}}
  
 \newcommand{\Cyc}{\mathfrak{C}} \newcommand{\Dih}{\mathfrak{D}}
 \newcommand{\Quat}{\mathfrak{Q}}

 \DeclareMathAlphabet{\cat}{OT1}{cmss}{m}{sl}
 \DeclareMathOperator{\tens}{\otimes}

 \DeclareMathOperator{\Gal}{Gal} \newcommand{\Set}{\cat{Set}}
 \newcommand{\Cov}{\cat{Cov}} \newcommand{\Et}{\cat{\acute Et}}
  \newcommand{\Etex}{\cat{\acute Etex}}

  \newcommand{\Sym}{\mathfrak{S}}
 \newcommand{\A}{\mathfrak{A}} 
 \DeclareMathOperator{\Id}{Id}
 \DeclareMathOperator{\funcX}{\mathbf{X}}
 \DeclareMathOperator{\Inv}{Inv} \DeclareMathOperator{\Int}{Int}
 \DeclareMathOperator{\sgn}{sgn} \DeclareMathOperator{\diag}{Diag}
 \DeclareMathOperator{\Iso}{Iso}
 
 \DeclareMathOperator{\id}{Id} 
 
 \DeclareMathOperator{\Hom}{Hom} 
 \DeclareMathOperator{\disc}{Disc}
 \DeclareMathOperator{\Res}{Res}
  \DeclareMathOperator{\Skew}{Skew}
 \DeclareMathOperator{\SSym}{Sym} 
  \DeclareMathOperator{\Fix}{Fix}
 \DeclareMathOperator{\Aut}{Aut} 
  
  \DeclareMathOperator{\Tr}{Tr}
 \DeclareMathOperator{\Spin}{Spin} 
  
  \DeclareMathOperator{\PGO}{PGO}
 \DeclareMathOperator{\OO}{O}

 \DeclareMathOperator{\funcM}{\mathbf{M}}

\newcommand{\joinrelshort}{\mathrel{\mkern-9mu}}
\newcommand{\shortlongrightarrow}{\relbar\joinrelshort\rightarrow}
\newcommand{\iso}{\mathrel{\mathop{\setbox0\hbox{$\mathsurround0pt
        \shortlongrightarrow$}\ht0=0.7\ht0\box0}\limits
    ^{\sim\mkern2mu}}}
   \newcommand{\ljoinrelshort}{\mathrel{\mkern-16mu}}
\newcommand{\lshortlongrightarrow}{\relbar\ljoinrelshort\leftarrow}
\newcommand{\liso}{\mathrel{\mathop{\setbox0\hbox{$\mathsurround0pt
        \lshortlongrightarrow$}\ht0=0.7\ht0\box0}\limits
    ^{\sim\mkern2mu}}}

\newcommand{\hcenter}[1]{{\leavevmode\setbox0\hbox{#1}\kern-.5\wd0\box0}}

\DeclareMathAlphabet{\cat}{OT1}{cmss}{m}{sly}

\theoremstyle{plain} %
\numberwithin{equation}{section}

\newtheorem{thm}[equation]{Theorem}

\newtheorem{prop}[equation]{Proposition}
\newtheorem{cor}[equation]{Corollary}

\theoremstyle{definition} %
\numberwithin{equation}{section}
\newtheorem{example}[equation]{Example}

\newtheorem{remark}[equation]{Remark}

\newcommand{\card}[1]{\lvert#1\rvert}

\newcounter{eqalignnumcnt}

\newcommand{\OLDref}{} \let\OLDref\ref
\renewcommand{\ref}[1]{\textup{\OLDref{#1}}}

\newcommand{\nop}[3]{}

\dedicatory{Dedicated with great friendship to R. Parimala at the
  occasion of her $60^{th}$ birthday}

\title{Triality and \'etale algebras} \author{Max-Albert Knus \and
  Jean-Pierre Tignol}
\address{Department Mathematik\\
  ETH Zentrum\\
  CH-8092 Z\"urich\\
  Switzerland} \email{knus@math.ethz.ch}
\address{Institut de Math\'ematique Pure et Appliqu\'ee\\
  Universit\'e catholique de Louvain\\
  B-1348 Louvain-la-Neuve\\
  Belgium} \email{jean-pierre.tignol@uclouvain.be} \thanks{The second
  author is supported in part by the F.R.S.--FNRS (Belgium)}

\begin{document}

\begin{abstract}
  Trialitarian automorphisms are related to automorphisms of order~$3$
  of the Dynkin diagram of type $D_4$. Octic \'etale algebras with
  trivial discriminant, containing quartic subalgebras, are classified
  by Galois cohomology with value in the Weyl group of type $D_4$.
  This paper discusses triality for such \'etale extensions.
\end{abstract}
\maketitle
\section{Introduction}
All Dynkin diagrams but one admit at most automorphisms of
order two, which are related to duality in algebra and geometry.
The Dynkin diagram \-
of $D_4$\\[1ex]
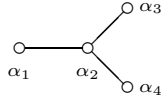
\[
\unitlength =0.5ex
\begin{picture}(-10,0)%
        \put(11,7){$\vcenter{\hbox{$\scriptstyle \alpha _3$}}$}%
        \put(11,-7){$\vcenter{\hbox{$\scriptstyle \alpha_ 4
\phantom{\scriptstyle -1}$}}$}%
\end{picture}%
\vcenter{\hbox{\begin{picture}(0,0)%
      \put(0,0){\circle{2}}%
      \put(0,-5){\hcenter{$\scriptstyle \alpha _1$}}
        \put(1,0){\line(1,0){10}}%
        \put(12,0){\circle{2}}%
        \put(12,-5){\hcenter{$\scriptstyle \alpha_ 2$}}%
                \put(19,7){\circle{2}}%
        \put(19,-7){\circle{2}}%
      \put(12.71,0.71){\line(1,1){5.58}}
        \put(12.71,-0.71){\line(1,-1){5.58}}
              \end{picture}}}
\]\\[0.1ex]

is special, in the sense that it admits automorphisms of order $3$.
Algebraic and geometric objects related to $D_4$ are of particular
interest as they also usually admit exceptional automorphisms of order
$3$, which are called trialitarian. For example the special projective
orthogonal group $\PGO^+_8$ or the simply connected group $\Spin_8$
admit outer automorphisms of order $3$. As already observed by E.
Cartan,~\cite{cartan}, the Weyl group $W(D_4) = \Sym_2^3 \rtimes
\Sym_4$ of $\Spin_8$ or of $\PGO^+_8$ similarly admits trialitarian
automorphisms.  Let $F$ be a field and let $F_s$ be a separable
closure of $F$. The Galois cohomology set $H^1\big(\Gamma,
W(D_4)\big)$, where $\Gamma$ is the absolute Galois group
$\Gal(F_s/F)$, classifies isomorphism classes of \'etale extensions 
$S/S_0$ where $S$ has dimension $8$, $S_0$ dimension $4$ 
and $S$ has trivial discriminant (see \S\ref{sec:cohomology}). 
There is an induced trialitarian action on  $H^1\big(\Gamma,
W(D_4)\big)$, which associates to the isomorphism class of an extension $S/S_0$ as above, two extensions
$S'/S_0'$ and $S''/S''_0$, of the same kind, so that the triple
$(S/S_0,S'/S'_0,S''/S''_0)$ is cyclically permuted by triality.
 This paper is devoted to
the study of such triples of \'etale algebras. It grew out of a study,
in the spirit of \cite{Tits57}, of Severi-Brauer varieties over the 
``field of one element'', \cite{KT09}, which is in preparation (see also
\cite{Tits59} and \cite{Tits74}).

In Part~2 we describe some basic constructions on finite $\text{$\Gamma$-sets}$ and
\'etale algebras. Some results are well-known, others
were taken from \cite{KT09}, like the Clifford construction. In Section~3
we recall how $\Gamma$-sets and \'etale algebras are related to
 Galois cohomology. Section~4 is devoted to triality in connection with
$\Gamma$-sets and in Section~5 we discuss trialitarian automorphisms of
the  Weyl group $W(D_4)$. In Section~6  we consider triality at the level of
\'etale algebras.  We give in Table 1 a list of  isomorphism classes of \'etale algebras 
corresponding to the conjugacy
classes of subgroups of $W(D_4)$, together with a description of the triality
action. We also consider  \'etale algebras associated to subgroups of
 $W(D_4)$ which are fixed under triality. We then view in Section 7 triality as a
 way to create resolvents and give explicit formulae for polynomials defining
 \'etale algebras.
 Finally we give in the last section
 results of Serre on Witt invariants of $W(D_4)$.

We are grateful to Parimala for her unshakable interest in triality,
in particular for many discussions at earlier stages of this work and
we specially thank J-P.~Serre for communicating to us his results on Witt
and cohomological invariants of the group $W(D_4)$.  We also thank
Emmanuel Kowalski who introduced us to Magma \cite{magma} with much
patience, Jean Barge for his help with Galois cohomology and J. E.
Humphreys and B. M\"uhlherr for the reference to the paper
\cite{franzsen_howlett}. The paper \cite{JonesRoberts} on octic fields
was a very useful source of inspiration.
 Finally we are highly thankful to the referee
for many improvements.

\section{\'Etale algebras and $\Gamma$-sets}
\label{subsec:etale}
Throughout most of this work, $F$ is an arbitrary field. We denote by
$F_s$ a separable closure of $F$ and by $\Gamma$ the absolute Galois
group $\Gamma=\Gal(F_s/F)$, which is a profinite group.

A finite-dimensional commutative $F$-algebra $S$ is called
\emph{\'etale} (over $F$) if $S\otimes_FF_s$ is isomorphic to the
$F_s$-algebra $F_s^n=F_s\times\cdots\times F_s$ ($n$ factors) for some
$n\geq 1$.  \'Etale $F$-algebras are the direct products of finite separable field
extensions of $F$. We refer to \cite[\S18.A]{KMRT} for various equivalent
characterizations. \'Etale algebras (with $F$-algebra homomorphisms)
form a category $\Et_F$ in which finite direct products and finite
direct sums (= tensor products) are defined.
 
Finite sets with a continuous left action of $\Gamma$ (for the
discrete topology) are called (finite) $\Gamma$-sets. They form a
category $\Set_\Gamma$ whose morphisms are the $\Gamma$-equivariant
maps. Finite direct products and direct sums (= disjoint unions) are
defined in this category. We denote by $\lvert X\rvert$ the
cardinality of any finite set $X$.

For any \'etale $F$-algebra $S$ of dimension $n$, the set of
$F$-algebra homomorphisms
\[
\funcX(S)=\Hom_\text{$F$-alg}(S, F_s)
\]
is a $\Gamma$-set of $n$ elements since $\Gamma$ acts on $F_s$.
Conversely, if $X$ is a $\Gamma$-set of $n$ elements, the $F$-algebra
$\funcM(X)$ of $\Gamma$-equivariant maps $X\to F_s$ is an \'etale
$F$-algebra of dimension~$n$,
\[
\funcM(X)=\{f\colon X\to F_s\mid \gamma\bigl(f(x)\bigr)= f({}^\gamma
x) \text{ for $\gamma\in\Gamma$, $x\in X$}\}.
\]
As first observed by Grothendieck, there are canonical isomorphisms
\[
\funcM\bigl(\funcX(S)\bigr)\cong S,\qquad
\funcX\bigl(\funcM(X)\bigr)\cong X,
\]
so that the functors $\funcM$ and $\funcX$ define an anti-equivalence
of categories
\begin{equation}
  \label{eq:anteq}
  \Set_\Gamma\equiv\Et_F
\end{equation}
(see \cite[Proposition (4.3), p.~25]{deligne} or \cite[(18.4)]{KMRT}).
Under this anti-equivalence, the cardinality of $\Gamma$-sets
corresponds to the dimension of \'etale $F$-algebras, the disjoint
union $\sqcup$ in $\Set_\Gamma$ corresponds to the product $\times $
in $\Et_F$, and the product $\times$ in $\Set_\Gamma$ to the tensor
product $\tens$ in $\Et_F$. For any integer $n\geq1$, we let $\Et^n_F$
denote the groupoid\footnote{A groupoid is a category in which all 
  morphisms are isomorphisms.} whose objects are $n$-dimensional
\'etale $F$-algebras and whose morphisms are $F$-algebra isomorphisms,
and $\Set^n_\Gamma$ the groupoid of $\Gamma$-sets with $n$ elements.
The anti-equivalence \eqref{eq:anteq} restricts to an anti-equivalence
$\Set_\Gamma^n \equiv\Et_F^n$. The split \'etale algebra $F^n$
corresponds to the $\Gamma$-set $\boldsymbol{n}$ of $n$ elements with
trivial $\Gamma$-action.  \'Etale algebras of dimension~$2$ are also
called \emph{quadratic} \'etale algebras.

A morphism\footnote{We let morphisms of $\Gamma$-sets act on the right
  of the arguments (with the exponential notation) and use the usual
  function notation for morphisms in the anti-equivalent category of
  \'etale algebras.} of $\Gamma$-sets $Y\xleftarrow{\pi} Z$ is called a
\emph{$\Gamma$-covering} if the number of elements in each fiber
$y^{\pi^{-1}}\subset Z$ does not depend on $y\in Y$. This number is
called the \emph{degree} of the covering. For $n$, $d\geq1$ we let
$\Cov_\Gamma^{d/n}$ denote the groupoid whose objects are coverings of
degree $d$ of a $\Gamma$-set of $n$ elements and whose morphisms are
isomorphisms of $\Gamma$-coverings.

A homomorphism $S\xrightarrow{\varepsilon}T$ of \'etale $F$-algebras
is said to be an \emph{extension of degree $d$ of \'etale algebras} if
$\varepsilon$ endows $T$ with a structure of a free $S$-module of
rank~$d$.  This corresponds under the
anti-equivalence~\eqref{eq:anteq} to a \emph{covering of degree~$d$}:
\[
\funcX(S)\xleftarrow{\funcX(\varepsilon)} \funcX(T)
\]
(see \cite{KT09}). Let $\Etex_F^{d/n}$ denote the groupoid of \'etale
extensions $S\xrightarrow{\varepsilon}T$ of degree~$d$ of $F$-algebras
with $\dim_FS=n$ (hence $\dim_FT=dn$). \relax From~\eqref{eq:anteq} we
obtain an anti-equivalence of groupoids
\[
\Etex^{d/n}_F\equiv\Cov^{d/n}_\Gamma.
\]

The $\Gamma$-covering with trivial $\Gamma$-action
\begin{equation}\label{equ:trivialcovering}
\boldsymbol{d}/\boldsymbol{n}:\quad \boldsymbol{n}\xleftarrow{p_1}
\boldsymbol{n}\times\boldsymbol{d}
\end{equation}
where $p_1$ is the first projection corresponds to the extension
$F^n \to (F^d)^n$.

Of particular importance in the sequel are coverings of degree~$2$,
which are also called \emph{double coverings}. Each such covering
$Y\xleftarrow{\pi}Z$ defines a canonical automorphism 
$Z\xleftarrow{\sigma} Z$ of order~$2$, which interchanges the elements in each fiber of
$\pi$. Clearly, this automorphism has no fixed points. Conversely, if
$Z$ is any $\Gamma$-set and $ Z\xleftarrow{\sigma} Z$ is an automorphism
of order~$2$ without fixed points, the set of orbits
\[
Z/\sigma=\bigl\{\{z,z^\sigma\}\mid z\in Z\}
\]
is a $\Gamma$-set and the canonical map $(Z/\sigma)\leftarrow Z$ is a
double covering.  An \emph{involution} of a $\Gamma$-set with an
even number of elements is any automorphism of order~$2$ without fixed
points.

Let $\sigma\colon S\to S$ be an automorphism of order~$2$ of an
\'etale $F$-algebra $S$, and let $S^\sigma\subset S$ denote the
$F$-subalgebra of fixed elements, which is necessarily \'etale. The
following conditions are equivalent (see \cite{KT09}):
  \begin{enumerate}
  \item[(a)] the inclusion $S^\sigma\to S$ is a quadratic \'etale
    extension of $F$-algebras;
  \item[(b)] the automorphism $\funcX(\sigma)$ is an involution on
    $\funcX(S)$.
  \end{enumerate}
  We say under these equivalent conditions that the automorphism
  $\sigma$ is an \emph{involution} of the \'etale $F$-algebra $S$.

\subsection{Basic constructions on $\Gamma$-sets}
\label{subsub:discset}
We recall from \cite[\S18]{KMRT} and \cite[\S2.1]{knustignol} the
construction of the discriminant $\Delta(X)$ of a $\Gamma$-set $X$
with $\lvert X\rvert=n\geq2$. Consider the set of $n$-tuples of 
elements in $X$:
\[
\Sigma_n(X)=\{(x_1,\ldots, x_n)\mid X=\{x_1,\ldots, x_n\}\}.
\]
This $\Gamma$-set carries an obvious transitive (right) action of the
symmetric group $\Sym_n$. The \emph{discriminant} $\Delta(X)$ is the
set of orbits under the alternating group $\A_n$:
\[
\Delta(X)=\Sigma_n(X)/\A_n.
\]
It is a $\Gamma$-set of two elements, so $\Delta$ is a functor
\[
\Delta\colon\Set^n_\Gamma\to\Set^2_\Gamma.
\]

For any covering $Z_0\xleftarrow{\pi}Z$ of degree~$2$ with $\lvert Z_0
\rvert=n$, (hence $\card{Z}=2n$), we consider the set of (not
necessarily $\Gamma$-equivariant) sections of~$\pi$:
\[
C(Z/Z_0 )=\bigl\{\{z_1,\ldots,z_n\}\subset Z\mid \{z_1^\pi,\ldots,
z_n^\pi\}=Z_0 \bigr\}.
\]
It is a $\Gamma$-set with $2^n$ elements, so $C$ is a functor
\[
C\colon\Cov_\Gamma^{2/n}\to\Set_\Gamma^{2^n},
\]
called the \emph{Clifford functor} (see~\cite{KT09}). The $\Gamma$-set
$C(Z/Z_0 )$ is equipped with a canonical surjective morphism
\begin{equation}
\label{eq:delta}
\Delta(Z)\xleftarrow{\delta} C(Z/Z_0),
\end{equation}
which is defined in \cite[\S2.2]{knustignol} as follows: let $\sigma
\colon Z\to Z$ be the involution canonically associated to the double
covering $Z_0 \xleftarrow{\pi}Z$, so the fiber of $z^\pi$ is
$\{z,z^{\sigma}\}$ for each $z\in Z$; then $\delta$ maps each section
$\{z_1,\ldots, z_n\}$ to the $\A_{2n}$-orbit of the $2n$-tuple
$(z_1,\ldots, z_n,z_1^{\sigma},\ldots, z_n^{\sigma})$,
\[
\{z_1,\ldots, z_n\}^\delta= (z_1,\ldots, z_n,z_1^{\sigma},\ldots,
z_n^{\sigma})^{\A_{2n}}.
\]
Note that the canonical involution $\sigma$ induces an involution
$\underline{\sigma}$ on $C(Z/Z_0)$, which maps each section $\omega$
to its complement $Z\setminus\omega$. We may view $C(Z/Z_0)$ as a
covering of degree~$2$ of the set of orbits
$C(Z/Z_0)/\underline{\sigma}$, and thus consider the Clifford
construction as a functor
\begin{equation} \label{def:definitionC}
C\colon\Cov_\Gamma^{2/n}\to\Cov_\Gamma^{2/2^{n-1}}.
\end{equation}

\begin{prop}
  \label{prop:delta}
  For sections $\omega$, $\omega'\in C(Z/Z_0 )$, we have
  $\omega^\delta=(\omega')^\delta$ if and only if
  $\card{\omega\cap\omega'}\equiv n\bmod2$. Moreover, denoting by
  $\iota$ the nontrivial automorphism of $\Delta(Z)$, we have
  \[
  \underline{\sigma} \circ\delta=
  \begin{cases}
    \delta&\text{if $n$ is even},\\
    \delta\circ\iota&\text{if $n$ is odd}.
  \end{cases}
  \]
\end{prop}

\begin{proof}
  Let $\omega=\{z_1,\ldots, z_n\}$ and $\omega'=\{z_1, \ldots, z_r,
  z_{r+1}^{{\sigma}}, \ldots, z_n^{\sigma}\}$, so
  $r=\card{\omega\cap\omega'}$,
  \begin{align*}
    \omega^\delta & =(z_1,\ldots, z_n,z_1^{\sigma},\ldots,
    z_n^{\sigma})^{\A_{2n}}\\
    \intertext{and}
    (\omega')^\delta & =(z_1, \ldots, z_r, z_{r+1}^{\sigma}, \ldots,
    z_n^{\sigma}, z_1^{\sigma}, \ldots, z_r^{\sigma}, z_{r+1}, \ldots,
    z_n)^{\A_{2n}}.
  \end{align*}
  The permutation $\sigma'$ that interchanges $z_i$ and $z_i^{\sigma}$
  for $i=r+1$, \ldots, $n$ satisfies
  \[
  (z_1, \ldots, z_n, z_1^{\sigma}, \ldots, z_n^{\sigma})^{\sigma'}=
  (z_1, \ldots, z_r, z_{r+1}^{\sigma}, \ldots, z_n^{\sigma},
  z_1^{\sigma}, \ldots, z_r^{\sigma}, z_{r+1}, \ldots, z_n);
  \]
  it is in $\A_{2n}$ if and only if $n-r$ is even, which means
  $\card{\omega\cap\omega'}\equiv n\bmod2$.  For
  $\omega'=\omega^{\underline{\sigma}}$ the complement of $\omega$ we
  have $\card{\omega\cap\omega^{\underline{\sigma}}}=0$, hence
  $\omega^{\underline{\sigma}\delta}=\omega^\delta$ if and only if
  $n\equiv 0\bmod2$.
\end{proof}

\subsection{Oriented  $\Gamma$-sets}
An \emph{oriented $\Gamma$-set} is a pair $(Z, \partial_Z)$ where $Z$
is a $\Gamma$-set and $\partial_Z$ is a fixed isomorphism of
$\Gamma$-sets $ \boldsymbol{2} \liso \Delta(Z) $. In particular the
$\Gamma$-action on $\Delta(Z)$ is trivial. There are two possible
choices for $\partial_Z$. A choice is an \emph{orientation of~$Z$}.
Oriented $\Gamma$-sets with $n$ elements form a groupoid
$(\Set_{\Gamma}^n)^+$ whose morphisms are isomorphisms $Z_2
\xleftarrow{f} Z_1$ such that $  \Delta(f) \circ \partial_{Z_2} = \partial_{Z_1}$.
Similarly \emph{oriented coverings} are pairs $(Z/Z_0,
\partial_Z)$ where $Z_0 \leftarrow Z$ is a $\Gamma$-covering and
$\partial_Z$ is an orientation of~$Z$. We denote by
$(\Cov_\Gamma^{d/n})^+$ the groupoid of oriented coverings of degree
$d$ of $\Gamma$-sets with $n$ elements. Changing the orientation
through the twist $\boldsymbol{2} \xleftarrow{\iota} \boldsymbol{2}$
defines an involutive functor
 \[
 \kappa\colon (\Cov_\Gamma^{d/n})^+ \to (\Cov_\Gamma^{d/n})^+.
\]
\begin{prop}\label{prop:Cinpair}
  If $n$ is even the functor $C\colon \Cov_\Gamma^{2/n} \to
  \Cov_\Gamma^{2/2^{n-1}}$ of \eqref{def:definitionC} restricts to a
  pair of functors
\[
C_1,\,C_2 \colon (\Cov_\Gamma^{2/n})^+ \to \Cov_\Gamma^{2/2^{n-2} }.
\]
Moreover two sections $\omega$ and $\omega'$ of the oriented $\Gamma$-covering 
$(Z/Z_0,\partial_Z)$
lie in the same set $C_1(Z/Z_0,\partial_Z )$ or $C_2(Z/Z_0,\partial_Z )$ if and only if
$\vert\omega \cap \omega' \vert \equiv 0 \mod 2$.
\end{prop}

\begin{proof}
  Let $Z/Z_0$ be a $2/n$-covering. Proposition~\ref{prop:delta}
  implies that the covering $\Delta(Z) \xleftarrow{\delta} C(Z/Z_0
  ) $ factors through $C(Z/Z_0 )/ \underline{\sigma}$, where
  $\underline{\sigma}$ is the canonical involution of $C(Z/Z_0 )$:
  \[
  \Delta(Z) \leftarrow C(Z/Z_0 )/\underline{\sigma} \leftarrow
  C(Z/Z_0).
 \]
 Thus, if $Z/Z_0$ is oriented, we may use the given isomorphism
 $ \boldsymbol{2} \xleftarrow{\partial_Z} \Delta(Z)$ to define the
 $\Gamma$-sets
 \begin{align*}
  C_1(Z/Z_0,\partial_Z) & =\{\omega\in C(Z/Z_0)\mid
  \omega^{\delta\partial_Z}=1\}\\
  \intertext{and}
  C_2(Z/Z_0,\partial_Z) & =\{\omega\in C(Z/Z_0)\mid
  \omega^{\delta\partial_Z}=2\}.
  \end{align*}
  Obviously, we have $
 C(Z/Z_0 ) = C_1(Z/Z_0,\partial_Z) \sqcup C_2(Z/Z_0,\partial_Z) 
 $, and Proposition~\ref{prop:delta} shows that $\underline{\sigma}$
 restricts to involutions on
 $C_1(Z/Z_0,\partial_Z)$ and $ C_2(Z/Z_0,\partial_Z)$. The last claim
 also follows from Proposition~\ref{prop:delta}.
  \end{proof}
  
  We call the two functors $C_1$ and $C_2$ the \emph{spinor
    functors}. Note that when $n$ is even an orientation $\partial_Z$
  on $Z/Z_0$ can also be defined by specifying whether a given section
  $\omega\in C(Z/Z_0)$ lies in $C_1(Z/Z_0,\partial_Z)$ or
  $C_2(Z/Z_0,\partial_Z)$. Indeed, $\omega\in C_1(Z/Z_0,\partial_Z)$
  if and only if $\omega^\delta\in\Delta(Z)$ is mapped to $1$, which
  determines $\partial_Z$ uniquely. We shall avail ourselves of this
  possibility to define orientations on coverings in
  $\Cov_\Gamma^{2/4}$ in \S\ref{sec:trialcov}.

\subsection{Basic constructions on \'etale algebras}
We now consider analogues of the functors $\Delta$ and $C$ for \'etale
algebras and \'etale extensions.

For $S$ an \'etale $F$-algebra of dimension $n\geq2$, the discriminant
$\Delta(S)$ is a quadratic \'etale $F$-algebra such that

\[
\funcX\bigl(\Delta(S)\bigr) = \Delta\bigl(\funcX(S)\bigr).
\]
We thus have a functor
\[
\Delta\colon\Et^n_F\to\Et^2_F\qquad\text{for $n\geq2$}.
\]
If the field $F$ has characteristic different from~$2$, it is usual to
represent $\Delta(S)$ as $F[x]/\big(x^2 -\disc(S)\big),\, \disc(S) \in
F^\times$, and the class of $\disc(S)$ in $F^\times/(F^{\times})^2$ is
the usual discriminant. We refer to \cite[p.~291--293]{KMRT} and
\cite[\S3.1]{knustignol} for details.

Let $S\xrightarrow{\varepsilon}T$ be an \'etale extension of degree
$2$ of (\'etale) $F$-algebras, with $\dim_FS=n$, $\dim_FT=2n$. In
\cite[\S3.2]{knustignol}\footnote{The notation $\Omega$ is used for $C$ in
  \cite{knustignol}.} we define an \'etale $F$-algebra $C(T/S)$ such
that
\[
\funcX\bigl(C(T/S)\bigr)=C\bigl(\funcX(T)/\funcX(S)\bigr).
\]
\begin{example} If $\dim_FT=2$ and $S=F$, we have $C(T/S)=T$.

 For $S_1$, $S_2$ \'etale
  algebras of arbitrary dimension, and for arbitrary \'etale extensions
  $T_1/S_1$ and $T_2/S_2$ of degree $2$, there is a canonical isomorphism
\[
P\colon C\bigl( (T_1 \times T_2) / (S_1 \times S_2)\bigr) \ \iso \ C(
T_1/S_1) \otimes C( T_2/S_2).
\]
\end{example}
We call the $2^n$-dimensional algebra $C(T/S)$ the \emph{Clifford
  algebra} of $T/S$. It admits a canonical involution
$\underline{\sigma}$. If $\dim_FS$ is even $\underline{\sigma}$ is the
identity on $\Delta(T)$. The canonical morphism $\delta$ of
\eqref{eq:delta}
\[
\Delta\bigl(\funcX(T)\bigr)\xleftarrow{\delta}
 C\bigl(\funcX(T)/\funcX(S)\bigr)
\]
yields a canonical $F$-algebra homomorphism which we again denote by
$\delta$,
\[
\Delta(T)\xrightarrow{\delta} C(T/S),
\]
so that $C(T/S)$ is an \'etale extension of degree~$2^{n-1}$ of a
quadratic \'etale $F$-algebra.

\subsection{Oriented \'etale algebras}
As for oriented $\Gamma$-sets we define \emph{oriented \'etale
  algebras} as pairs $(S,\partial_S)$ where $S$ is an \'etale algebra
and $\partial_S\colon \Delta(S) \iso F \times F$ is an isomorphism of
$F$-algebras. \emph{Oriented extensions of \'etale algebras} are pairs
$(S/S_0,\partial_S)$ where $S/S_0$ is an extension of \'etale algebras
and $\partial_S\colon \Delta(S) \iso F \times F$ is an isomorphism of
$F$-algebras.  We have corresponding groupoids $(\Et_F^n)^+$,
$(\Etex_F^{d/n})^+$ and anti-equivalences
 \[
 (\Set_\Gamma^n)^+ \equiv (\Et_F^n)^+\ \text{and} \ 
 (\Cov_\Gamma^{d/n})^+ \equiv (\Etex_F^{d/n})^+.
 \]
 Switching the orientation induces an involutive functor $\kappa$ on
 these groupoids.
 
 The Clifford functor $C$ restricts to a pair of \emph{spinor
   functors}
 \begin{equation}
  \label{eq:Cetaleplus}
 C_1, C_2 \colon (\Etex_F^{2/n})^+ \to  \Etex_F^{2/2^{n-2} }
 \end{equation}
 if $n$ is even.

\begin{remark} The terminology used above owes its origin to the fact that
  the Clifford functor is related to the theory of Clifford algebras
  in the framework of quadratic forms and central simple algebras with
  involution. We refer to \cite{KT09} for details and more properties
  of the Clifford construction.
\end{remark}

\section{Cohomology}
\label{sec:cohomology}

For any integer $n\geq1$, we consider the $\Gamma$-set
$\boldsymbol{n}=\{1,\ldots, n\}$ with the trivial $\Gamma$-action and
let $\Sym_n$ denote the symmetric group on $\boldsymbol{n}$, i.e., the
automorphism group of $\boldsymbol{n}$,
\[
\Sym_n=\Aut(\boldsymbol{n}).
\]
Recall from \cite[\S28.A]{KMRT} that the cohomology set
$H^1(\Gamma,\Sym_n)$ (for the trivial action of $\Gamma$ on $\Sym_n$)
is the set of continuous group homomorphisms $\Gamma\to\Sym_n$ 
(``cocycles'') up to
conjugation. 

Letting $\Iso(\Set^n_\Gamma)$ denote the set of isomorphism classes in
$\Set^n_\Gamma$, we have a canonical bijection of pointed sets
\begin{equation}
\label{eq:setH1}
\Iso(\Set^n_\Gamma)\iso H^1(\Gamma,\Sym_n).
\end{equation}

Cohomology sets can also be used to describe isomorphism classes of
$\Gamma$-coverings: for any integers $n$, $d\geq1$, the group of
automorphisms of the $\Gamma$-covering with trivial
$\Gamma$-action~$\boldsymbol{d}/\boldsymbol{n}$ is the wreath product
(of order~$(d!)^nn!$)
\[
\Aut(\boldsymbol{d}/\boldsymbol{n})
=\Sym_d\wr\Sym_n\quad(=\Sym_d^n\rtimes \Sym_n).
\]
The same construction as above yields a canonical bijection
\begin{equation}
\label{eq:covH1}
\Iso(\Cov^{d/n}_\Gamma)\iso H^1(\Gamma,\Sym_d\wr\Sym_n),
\end{equation}
where the $\Gamma$-action on $\Sym_d\wr\Sym_n$ is trivial; see
\cite[\S4.2]{knustignol}.  The automorphism group of the oriented
$\Gamma$-covering $(\boldsymbol{d}/\boldsymbol{n},
\partial_{\boldsymbol{n}\times \boldsymbol{d}})$ is the group
\[
(\Sym_d \wr \Sym_n)^+ = ( \Sym_d \wr \Sym_n) \cap \A_{dn}
\]
so that
\begin{equation}
\label{eq:covplusH1}
\Iso\big((\Cov^{d/n}_\Gamma)^+\big) \iso H^1\big(\Gamma,(\Sym_d\wr\Sym_n)^+\big).
\end{equation}
 We now
assume that $\Gamma$ is the absolute Galois group $\Gamma=\Gal(F_s/F)$
of a field $F$ and use the notation $H^1(F,\Sym_n)$ for
$H^1(\Gamma,\Sym_n)$.  The anti-equivalence $\Set_\Gamma^n
\equiv\Et_F^n$ and the bijection \eqref{eq:setH1} induce canonical
bijections 
\[
\Iso(\Et^n_F)\cong\Iso(\Set^n_\Gamma) \cong H^1(F,\Sym_n)
\]
The bijection $\Iso(\Et^n_F)\cong H^1(F,\Sym_n)$ may of course
also be defined directly since
\[
\Aut_\text{$F$-alg}(F^{n})\cong\Sym_n,
\]
see \cite[(29.9)]{KMRT}. Similarly, it follows from~\eqref{eq:covH1}, 
\eqref{eq:covplusH1},
and the anti-equivalence of groupoids
$\Etex^{d/n}_F\equiv\Cov^{d/n}_\Gamma$,    
$(\Etex^{d/n}_F)^+\equiv (\Cov^{d/n}_\Gamma)^+$, that we have canonical
bijections of pointed sets:
\begin{equation} \label{equ:isoetxcov}
\Iso(\Etex^{d/n}_\Gamma)\cong\Iso(\Cov^{d/n}_\Gamma)\cong
H^1(F,\Sym_d\wr\Sym_n)
\end{equation}
and
\begin{equation} 
\label{equ:etexcovplus}
\Iso\big((\Etex^{d/n}_\Gamma)^+\big)\cong\Iso\big((\Cov^{d/n}_\Gamma)^+\big)\cong
H^1\big(F,(\Sym_d\wr\Sym_n)^+\big).
\end{equation}

\begin{remark} \label{resolvents:rem}
Any group homomorphism $\varphi \colon G \to H$, where
$G$ and $H$ are automorphism groups of finite sets or of finite double coverings,
 induces a map on the level of cocycles $\varphi_* \colon (\gamma \colon\Gamma \to G)
 \mapsto (\varphi\circ \gamma\colon \Gamma \to H)$. Thus $\varphi$
 associates in   a ``canonical way''  an \'etale algebra 
(or an \'etale algebra with involution) $E_\varphi$, whose isomorphism class  belongs to $H^1(F,H)$, to an \'etale algebra $E$ (or an \'etale algebra 
 $E$ with involution), whose class belongs to $H^1(F,G)$.
  We say that
 the algebra $E_\varphi$ is a \emph{resolvent} of $E$. For example the discriminant
 $\Delta(E)$ is the resolvent of $E$ associated to the parity map $\Sym_n \to \Sym_2$.
 Other examples of resolvents will be discussed in relation with triality.
 \end{remark}

\section{Triality and $\Gamma$-coverings}
\label{sec:trialcov}

Recall the functor $C$, which associates to any double covering its set of
sections. For oriented $2/4$-coverings of $\Gamma$-sets, it leads to two
functors
\[
C_1,\;C_2 \colon (\Cov_\Gamma^{2/4})^+ \to \Cov_\Gamma^{2/4},
\]
see Proposition~\ref{prop:Cinpair}.
The functors $C_1$ and $C_2$ together with the functor $\kappa$, which
changes the orientation, give an explicit description of an
action of the group $\Sym_3$ on the pointed set
$\Iso\bigl((\Cov_\Gamma^{2/4})^+\bigr)$.

\begin{thm} \label{thm:Candtriality}
  The functors $C_1$, $C_2 \colon (\Cov_\Gamma^{2/4})^+ \to
  \Cov_\Gamma^{2/4}$ factor through the forgetful functor
  $\mathcal{F}\colon(\Cov_\Gamma^{2/4})^+\to \Cov_\Gamma^{2/4}$, i.e.,
  there are functors 
  \[
  C_1^+,\; C_2^+\colon (\Cov_\Gamma^{2/4})^+ \to
  (\Cov_\Gamma^{2/4})^+
  \]
  such that $\mathcal{F}\circ C_i^+ = C_i$ for
  $i=1$, $2$. These functors satisfy natural equivalences:
  \[
  (C_1^+)^3 = \id,\qquad (C_1^+)^2=C_2,\qquad C_1^+\kappa=\kappa C_2^+.
  \]
\end{thm}

\begin{proof}
  Let $(Z/Z_0,\partial)$ be an object in $(\Cov_\Gamma^{2/4})^+$ and
  let $\sigma$ denote the involution of $Z/Z_0$. Consider a real
  vector space $V$ with basis $(e_1,e_2,e_3,e_4)$. Fixing a bijection
  $\varphi$ between a section $\omega\in C_1(Z/Z_0,\partial)$ and
  $\{e_1, \ldots, e_4\}$, we identify $Z$ with a subset of $V$ by
  \[
  z\mapsto
  \begin{cases}
    z^\varphi&\text{if $z\in\omega$,}\\
   -z^{\sigma\varphi}&\text{if $z\notin\omega$.}
  \end{cases}
  \]
  Thus, $Z=\{\pm e_1,\pm e_2, \pm e_3, \pm e_4\}$ and $\sigma$ acts on
  $Z$ by mapping each element to its opposite. The action of $\Gamma$
  on $Z$ extends to a linear action on $V$ since it commutes with
  $\sigma$. We also identify $C(Z/Z_0)$ with a subset of $V$ by the
  map
  \[
  \omega'\mapsto{\textstyle\frac12}\sum_{z\in\omega'}z.
  \]
  The set $C(Z/Z_0)$ then consists of the following vectors and their
  opposite:
  \begin{align*}
    f_1 & ={\textstyle\frac12}(e_1+e_2+e_3+e_4),& g_1 &
    ={\textstyle\frac12}(e_1+e_2+e_3-e_4),\\
    f_2 & ={\textstyle\frac12}(e_1+e_2-e_3-e_4),& g_2 &
    ={\textstyle\frac12}(e_1+e_2-e_3+e_4),\\
    f_3 & ={\textstyle\frac12}(e_1-e_2+e_3-e_4),& g_3 &
    ={\textstyle\frac12}(e_1-e_2+e_3+e_4),\\
    f_4 & ={\textstyle\frac12}(-e_1+e_2+e_3-e_4),& g_4 &
    ={\textstyle\frac12}(e_1-e_2-e_3-e_4).
  \end{align*}
  Since $\omega\in C_1(Z/Z_0,\partial)$, we have (see Proposition~\ref{prop:Cinpair}
  \begin{align*}
    C_1(Z/Z_0,\partial) & =\{\pm f_1,\pm f_2,\pm f_3, \pm f_4\}\\
    \intertext{and}
    C_2(Z/Z_0,\partial) & =\{\pm g_1, \pm g_2, \pm g_3, \pm g_4\}.
  \end{align*}
  The canonical involutions on $C_1(Z/Z_0,\partial)$ and
  $C_2(Z/Z_0,\partial)$ map each vector to its opposite.
  Note that these identifications are independent of the choice of the
  section~$\omega$ in $C_1(Z/Z_0,\partial)$ and of the bijection
  $\varphi$. 

  Let $\partial_1$ be the orientation of $C_1(Z/Z_0,\partial)$ such
  that $C_1(C_1(Z/Z_0,\partial),\partial_1)$ contains the section
  $\{f_1,f_2,f_3,f_4\}$. Thus, by Proposition~\ref{prop:Cinpair},
  $C_1(C_1(Z/Z_0,\partial),\partial_1)$ consists of the following
  sections:
  \[
  \pm\{f_1,f_2,f_3,f_4\},\;\pm\{f_1,f_2,-f_3,-f_4\},\;
  \pm\{f_1,-f_2,f_3,-f_4\},\;\pm\{-f_1,f_2,f_3,-f_4\}.
  \]
  They are characterized by the property that for each $i=1$, \ldots,
  $4$ they contain a section $\pm f_j$ containing $e_i$ and a section
  $\pm f_k$ containing $-e_i$. Identifying these sections to vectors
  in $V$ as above, we obtain
  \begin{multline*}
  C_1(C_1(Z/Z_0,\partial),\partial_1)=
  \{\pm{\textstyle\frac12}(f_1+f_2+f_3+f_4),
  \pm{\textstyle\frac12}(f_1+f_2-f_3-f_4),\\
  \pm{\textstyle\frac12}(f_1-f_2+f_3-f_4),
  \pm{\textstyle\frac12}(-f_1+f_2+f_3-f_4)\}.
  \end{multline*}
  Likewise, let $\partial_2$ be the orientation of
  $C_2(Z/Z_0,\partial)$ such that
  \begin{multline*}
  C_1(C_2(Z/Z_0,\partial),\partial_2)=
  \{\pm{\textstyle\frac12}(g_1+g_2+g_3+g_4),
  \pm{\textstyle\frac12}(g_1+g_2-g_3-g_4),\\
  \pm{\textstyle\frac12}(g_1-g_2+g_3-g_4),
  \pm{\textstyle\frac12}(-g_1+g_2+g_3-g_4)\}.
  \end{multline*}
  Define the functors $C_1^+$, $C_2^+\colon (\Cov_\Gamma^{2/4})^+ \to
  (\Cov_\Gamma^{2/4})^+$ by
  \[
  C_1^+(Z/Z_0,\partial) = (C_1(Z/Z_0),\partial_1)
  \qquad\text{and}\qquad C_2^+(Z/Z_0,\partial) =
  (C_2(Z/Z_0),\partial_2).
  \]
  By definition, it is clear that $\mathcal{F}\circ C_i^+=C_i$ for
  $i=1$, $2$. To establish the natural equivalences, consider the
  linear map $\mu\colon V\to V$ defined by $\mu(e_i)=f_i$ for $i=1$,
  \ldots, $4$. Using this map, we may rephrase the definition of
  $C_1^+$ as follows: for $Z=\{\pm e_1,\pm e_2, \pm e_3, \pm e_4\}$
  with the orientation $\partial$ such that $C_1(Z/Z_0,\partial)\ni
  \mu(e_1)$, we have 
  \[
  C_1^+(Z/Z_0,\partial) =\{\pm\mu(e_1),\pm\mu(e_2), \pm\mu(e_3),
  \pm\mu(e_4)\}
  \]
  with the orientation such that
  \[
  C_1\bigl( C_1^+(Z/Z_0,\partial)\bigr)\ni {\textstyle\frac12}\bigl(
  \mu(e_1)+\mu(e_2)+\mu(e_3)+\mu(e_4)\bigr).
  \]
  Note that ${\textstyle\frac12}\bigl(
  \mu(e_1)+\mu(e_2)+\mu(e_3)+\mu(e_4)\bigr)=
  \mu(f_1)=\mu^2(e_1)$. Therefore, substituting $\mu(e_i)$ for $e_i$,
  for $i=1$, \ldots, $4$, we obtain
  \[
  (C_1^+)^2(Z/Z_0,\partial) = \{\pm\mu^2(e_1), \pm\mu^2(e_2),
  \pm\mu^2(e_3), \pm\mu^2(e_4)\},
  \]
  endowed with an orientation such that
  \[
  C_1\bigl((C_1^+)^2(Z/Z_0,\partial)\bigr) \ni \mu^3(e_1).
  \]
  Computation shows that $\mu^2(e_i)=g_i$ for $i=1$, \ldots, $4$, and
  $\mu^3=\Id$. Since $\frac12(g_1+g_2+g_3+g_4)=e_1$, it follows that
  $(C_1^+)^2(Z/Z_0,\partial)= C_2^+(Z/Z_0,\partial)$. Similarly, we
  have
  \[
  (C_1^+)^3(Z/Z_0,\partial)= \{\pm\mu^3(e_1), \pm\mu^3(e_2), \pm\mu^3(e_3),
  \pm\mu^3(e_4)\} = Z,
  \]
  endowed with an orientation such that
  \[
  C_1\bigl((C_1^+)^3(Z/Z_0,\partial)\bigr)\ni\mu^4(e_1)=f_1,
  \]
  hence $(C_1^+)^3(Z/Z_0,\partial)=(Z/Z_0,\partial)$. Finally, we have
  \[
  \kappa(Z/Z_0,\partial) = \{\pm e_1, \pm e_2, \pm e_3, \pm e_4\}
  \]
  with an orientation such that $C_1\kappa(Z/Z_0,\partial)\ni
  \mu^2(e_1)$, hence
  \[
  C_1^+\kappa(Z/Z_0,\partial) = \{\pm\mu^2(e_1), \pm\mu^2(e_2),
  \pm\mu^2(e_3), \pm\mu^2(e_4)\}
  \]
  endowed with an orientation such that
  \[
  C_1\bigl(C_1^+\kappa(Z/Z_0,\partial)\bigr)\ni \mu^4(e_1)=f_1.
  \]
  Therefore, $C_1^+\kappa(Z/Z_0,\partial)=\kappa C_2^+(Z/Z_0,\partial)$.
\end{proof}

\begin{remark} \label{rem:hypercube}
The decomposition $C(Z/Z_0)=C_1(Z/Z_0,\partial)\sqcup
C_2(Z/Z_0,\partial)$ can also be viewed geometrically on a hypercube:
suppose $V=\mathbb{R}^4$ and let $(e_1,e_2,e_3,e_4)$ be the standard
basis. The set
\[
C(Z/Z_0)=\{{\textstyle\frac12}(\pm e_1\pm e_2\pm e_3\pm e_4)\}
\]
is the set of vertices of a hypercube $\mathcal{K}$ (see
Figure~1), and the set
\[
Z=\{\pm e_1,\pm e_2, \pm e_3, \pm e_4\}
\]
is in bijection with the set of $3$-dimensional cells of $\mathcal
K$.

\begin{center}    

{\scriptsize
\[
\unitlength =0.3mm
\begin{picture}(170,130)(5,5)
      \put(160,20){\circle{2}}
      \put(0,0){\circle{2}}
       \put(0,0){\line(1,0){100}}
        \put(0,0){\line(0,1){100}} 
        \put(0,0){\line(3,1){60}}
        \put(0,0){\line(4,3){39}}
         \put(-6,-7){$\vcenter{\hbox{$\mathsf  2$}}$}   
        \put(40,30){\line(1,0){50}}
         \put(40,30){\line(3,1){30}}
       \put(40,30){\circle{2}}
            \put(40,30){\line(0,1){50} }
                     \put(33,34){$\vcenter{\hbox{$\overline{\mathsf  5}$}}$}   
  \put(100,0){\line(-1,3){10}} 
  \put(100,0){\line(3,1){60}} 
       \put(100,0){\circle{2}}
        \put(100,0){\line(0,1){100}} 
         \put(104,-7){$\vcenter{\hbox{$\mathsf  3$}}$}   
       \put(60,20){\circle{2}}
        \put(60,20){\line(0,1){99}}
      \put(60,20){\line(1,0){99}}
               \put(63,13){$\vcenter{\hbox{$\mathsf  7$}}$}   
       \put(160,20){\circle{2}}
 \put(160,20){\line(0,1){100}}
    \put(160,20){\line(-2,1){41}}             
                 \put(162,13){$\vcenter{\hbox{$\mathsf  6$}}$}   
         \put(70,40){\circle{2}}
 \put(70,40){\line(0,1){50}}
  \put(70,40){\line(1,0){50}}
   \put(73,44){$\vcenter{\hbox{$ \overline{\mathsf 4}$}}$}  
       \put(90,30){\circle{2}}
        \put(90,30){\line(3,1){31}}
        \put(90,30){\line(0,1){50}}
         \put(94,25){$\vcenter{\hbox{$ \overline{\mathsf 8}$}}$}  
             \put(120,40){\circle{2}} 
          \put(120,40){\line(0,1){50}}
         \put(124,43){$\vcenter{\hbox{$ \overline{\mathsf 1}$}}$}   
     \put(- 6,105){$\vcenter{\hbox{$\mathsf  1$}}$}    
   \put(0,100){\line(2,-1){41}}%
  \put(0,100){\line( 100,0){100}}%
       \put(0,100){\circle{2}}
        \put(0,100){\line(3,1){60}}%
\put(60,120){\line(1,-3){10}}%
  \put(60,120){\line( 100,0){100}}%
       \put(60,120){\circle{2}}
        \put(63,125){$\vcenter{\hbox{$\mathsf  8$}}$}   
\put(40,80){\line(1,0){50}}%
  \put(40,80){\line( 3,1){30}}%
       \put(40,80){\circle{2}}
        \put(34,86){$\vcenter{\hbox{$ \overline{\mathsf 6}$}}$}  
       \put(70,90){\line(1,0){50}}%
       \put(70,90){\circle{2}}
 \put(73,83){$\vcenter{\hbox{$ \overline{\mathsf 3}$}}$}  
 \put(100,100){\line(3,1){59}}%
         \put(100,100){\circle{2}}
            \put(102,104){$\vcenter{\hbox{$\mathsf  4$}}$} 
      \put(90,80){\line(3,1){30}}%
       \put(90,80){\line(1,2){10}}%
               \put(90,80){\circle{2}}
                \put(93,75){$\vcenter{\hbox{$ \overline{\mathsf 7}$}}$}  
                      \put(120,90){\circle{2}}
\put(120,90){\line(4,3){40}}
 \put(124,84){$\vcenter{\hbox{$ \overline{\mathsf 2}$}}$}  
\put(160,120){\circle{2}}
               \put(162,125){$\vcenter{\hbox{$\mathsf  5$}}$}   
         \put(75,12){$\vcenter{\hbox{$ \overline{\mathsf C}$}}$}  
         \put(75,105){$\vcenter{\hbox{$ {\mathsf C}$}}$}  
                  \put(25,60){$\vcenter{\hbox{$ \overline{\mathsf B}$}}$}  
         \put(130,60){$\vcenter{\hbox{$ {\mathsf B}$}}$}  
                  \put(105,65){$\vcenter{\hbox{$ \overline{\mathsf A}$}}$}  
         \put(50,55){$\vcenter{\hbox{$ {\mathsf A}$}}$}  
                  \put(78,58){$\vcenter{\hbox{$ \overline{\mathsf D}$}}$}  
         \put(3,10){$\vcenter{\hbox{$ {\mathsf  D}$}}$}  
    \end{picture}    
   \]
   }
   
   \vspace{2ex}
   
  Figure 1
\end{center}

We identify
\[
Z =\{ A, \bar{A}, B, \bar{B},C, \bar{C},D, \bar{D} \}
\]
where $A, \ldots, \bar{C}$ are as in Figure~1, $D$ is the
big cell and $\bar{D}$ the small cell inside. The involution permutes a
cell with its opposite cell and the set $Z_0$ is obtained by
identifying pairs of opposite cells
\[
Z_0 =\{\{ A, \bar{A}\},\{ B, \bar{B}\}, \{C, \bar{C}\},\{D, \bar{D}\}
\}.
\]
To obtain a corresponding identification of $C(Z/Z_0)$ with the set of
vertices of $\mathcal K$, observe that a section of $Z/Z_0$ consists
of a set of four cells which are pairwise 
not opposite. Four such cells intersect in exactly one vertex
and conversely each vertex lies in four cells. With the
notation in Figure~1 we have the following identification:
\[
\begin{array}{llll}
1 =\{ A,\bar B,C,D\} & \bar{1} = \{\bar A,B,\bar C,\bar D\}& 
 2 =\{ A,\bar B,\bar C,D\} & \bar 2 =\{\bar A, B,C,\bar D\}\\
3 =\{ A, B,\bar C,D\} & \bar  {3} = \{\bar  A,\bar B, C,\bar D\} & 
 4 =\{ A, B,C,D\} &
 \bar 4 =\{\bar  A,\bar  B,\bar C,\bar  D\}\\
5 =\{\bar  A, B,C,D\} &  \bar {5} = \{ A,\bar B,\bar  C,\bar D\}& 
 6 =\{\bar  A, B,\bar C,D\} &
 \bar 6 =\{ A,\bar  B,C,\bar  D\}\\
7 =\{\bar  A, \bar B,\bar C,D\} &  \bar {7} = \{ A,B, C,\bar D\}& 
 8 =\{\bar  A,\bar  B, C,D\} &
 \bar 8 =\{ A, B, \bar C,\bar D\}\\
\end{array}
\]
This set of vertices decomposes into two
classes, two vertices being in the same class if the number of edges
in any path connecting them is even. One class is
\[ 
X = \{ 1,\bar 1,3,\bar 3,5,\bar5,7,\bar 7\}
\]
and the other
\[ 
Y = \{2,\bar 2,4,\bar 4,6,\bar 6,8,\bar 8 \}.
\]
We get coverings $X/X_0$ and $Y/Y_0$ by identifying opposite vertices
$v$ and $\bar v$. If $\Delta(Z)\simeq\boldsymbol{2}$, the
decomposition of $C(Z/Z_0)$ as the disjoint union $X/X_0 \sqcup Y/Y_0$
is $\Gamma$-compatible;
the functors $C_1$ and $C_2$ are given (up to a possible
permutation) by the rule
\[
C_1(Z/Z_0,\partial) = X/X_0 \qquad \text{and} \qquad C_2 (Z/Z_0,\partial)
= Y/Y_0.
\]
A section of $X/X_0$ is a set of four vertices in $X$ which are pairwise
not opposite. Four such vertices either lie on a
$3$-dimensional cell or are adjacent to exactly one vertex in the
complementary set $Y$.  A similar claim holds for a section of
$Y/Y_0$.  This leads to identifying:
\begin{equation} 
\begin{array}{lll} 
 A = \{ 1, 3,\bar  5,\bar  7\} = \{ 2, 4,\bar  6,\bar  8\} & 
 \bar  A = \{\bar  1,\bar  3, 5, 7\} = \{\bar  2,\bar  4, 6, 8\} \\
 B = \{\bar  1, 3, 5,\bar  7\} = \{ \bar  2, 4, 6,\bar  8\} &  
\bar  B = \{ 1,\bar  3, \bar 5, 7\} = \{ 2,\bar  4,\bar  6, 8\} \\
 C = \{ 1,\bar  3, 5, \bar 7\} = \{\bar  2, 4,\bar  6, 8\} & 
\bar  C = \{\bar  1, 3,\bar  5, 7\} = \{ 2,\bar  4, 6,\bar  8\} \\
 D = \{ 1, 3, 5, 7\} = \{ 2, 4, 6, 8\} &  
\bar  D = \{\bar  1,\bar  3,\bar  5,\bar  7\} = \{\bar  2,\bar  4,\bar
6,\bar  8\}  
\end{array} 
\end{equation}   
and
\begin{equation} 
\begin{array}{lll}

1 =\{ A,\bar B,C,D\}= \{ 2, 4,\bar 6, 8\} &  
\bar{1} = \{\bar A,B,\bar C,\bar D\} = \{\bar 2,\bar 4, 6,\bar 8\}\\ 
3 =\{ A, B,\bar C,D\}= \{ 2, 4, 6, \bar8\}  &  
\bar  {3} = \{\bar  A,\bar B, C,\bar D\}= \{\bar 2,\bar 4,\bar 6, 8\} \\ 
5 =\{\bar  A, B,C,D\}= \{\bar 2, 4, 6, 8\} &  
 \bar {5} = \{ A,\bar B,\bar  C,\bar D\}= \{ 2,\bar 4,\bar 6,\bar 8\}\\
7 =\{\bar  A, \bar B,\bar C,D\}= \{ 2,\bar 4, 6, 8\}  &  
 \bar {7} = \{ A,B, C,\bar D\}= \{\bar 2, 4,\bar 6,\bar 8\} \\[1ex] 
2 =\{ A,\bar B,\bar C,D\} = \{ 1, 3,\bar 5, 7  \} &  
\bar 2 =\{\bar A,B,C,\bar D\} = \{ \bar1,\bar 3, 5,\bar 7  \}\\ 
4 =\{ A, B, C,D\}= \{ 1, 3, 5,\bar 7  \} & 
 \bar 4 =\{\bar  A,\bar  B,\bar C,\bar  D\}= \{\bar 1,\bar 3,\bar 5, 7
 \} \\  
6 =\{\bar  A, B,\bar C,D\} = \{\bar 1, 3, 5, 7  \}&  
 \bar 6 =\{ A,\bar  B,C,\bar  D\} = \{ 1,\bar 3,\bar 5,\bar 7  \}\\ 
8 =\{\bar  A,\bar  B, C,D\}= \{ 1,\bar 3, 5, 7  \}  & 
 \bar 8 =\{ A, B, \bar C,\bar D\} =  \{\bar 1, 3,\bar 5, \bar7  \}, 
\end{array}
\end{equation}
hence the existence of decompositions $C(X/X_0) =Y/Y_0 \sqcup Z/Z_0$
and $C(Y/Y_0) =Z/Z_0 \sqcup X/X_0$ which, in fact, are decompositions
as $\Gamma$-sets.
\end{remark}

\begin{remark}
\label{rem:Hurwitz}
In the proof of Theorem~\ref{thm:Candtriality}, $\mu$ is not the unique linear map that can be
used to describe the $C_1^+$ and the $C_2^+$ construction. An
alternative description uses Hurwitz' quaternions. Choosing for $V$
the skew field of real quaternions $\mathbb H$ and for $(e_1, e_2,e_3,
e_4)$ the standard basis $(1,\,i,\,j,\,k)$, we have
\[
Z = \{ \pm 1,\, \pm i,\, \pm j,\, \pm k \},\qquad
C(Z/Z_0)=\{{\textstyle\frac12}(\pm1\pm i\pm j\pm k)\},
\]
so the union $Z \cup C(Z/Z_0) \subset \mathbb H$ is the
group $\mathbb H^1$ of Hurwitz integral quaternions of norm~$1$. The element
\[
\rho = -{\textstyle\frac{1}{2}}(1 +i+j+k)
\]
is of order $3$ in $\mathbb H^1$ and conjugation by $\rho$ permutes
$i$, $j$ and $k$ cyclically. The set $Z$ is in fact the underlying set
of the quaternionic group $\Quat_8$ and
\[
\mathbb H^1 = \Quat_8 \rtimes \Cyc_3
\]
where the cyclic group of three elements $\Cyc_3$ operates on $\Quat_8
$ via conjugation with $\rho$. If $\partial$ is the orientation of
$Z/Z_0$ such that $\rho\in C_1(Z/Z_0,\partial)$, we have
$C_1^+(Z/Z_0,\partial)=\rho\cdot Z$ with the orientation such that
$\rho^2\in C_1\bigl(C_1^+(Z/Z_0,\partial)\bigr)$, and
$C_2^+(Z/Z_0,\partial)=(C_1^+)^2(Z/Z_0,\partial)=\rho^2\cdot Z$ with
the orientation such that $1\in
C_1\bigl(C_2^+(Z/Z_0,\partial)\bigr)$. Note that, with respect to the
standard basis, multiplication by
$\rho$ is given by the matrix
\begin{equation}\label{matrixrho.equ}
\rho = \frac{1}{2}\left(
 \begin{matrix}
   -1 & 1 & 1 & 1 \\
   -1 & -1 & 1 & -1 \\
   -1 & -1 &- 1 & 1 \\
   -1 & 1 & -1 &-1
\end{matrix}\right)
\end{equation}
whereas the matrix of $\mu$ is
\begin{equation}\label{matrixmu.equ}
\mu = \frac{1}{2}\left( 
 \begin{matrix}
   1 & 1 & 1 & -1 \\
   1 & 1 & -1 & 1 \\
   1 & -1 & 1 & 1 \\
   1 & -1 & -1 & -1
\end{matrix}\right).
\end{equation}
\end{remark}

\section{The Weyl group of $D_4$}
\label{Weyl} 

The Dynkin diagram $D_4$
\begin{equation}\label{dynkin.equ}
\unitlength =0.5ex
\begin{picture}(-10,0)%
        \put(11,7){$\vcenter{\hbox{$\scriptstyle \alpha _3$}}$}%
        \put(11,-7){$\vcenter{\hbox{$\scriptstyle \alpha_ 4
\phantom{\scriptstyle -1}$}}$}%
\end{picture}%
\vcenter{\hbox{\begin{picture}(0,0)%
      \put(0,0){\circle{2}}%
      \put(0,-5){\hcenter{$\scriptstyle \alpha _1$}}
        \put(1,0){\line(1,0){10}}%
        \put(12,0){\circle{2}}%
        \put(12,-5){\hcenter{$\scriptstyle \alpha_ 2$}}%
                \put(19,7){\circle{2}}%
        \put(19,-7){\circle{2}}%
      \put(12.71,0.71){\line(1,1){5.58}}
        \put(12.71,-0.71){\line(1,-1){5.58}}
              \end{picture}}}
\end{equation}
\vspace{2ex}
             
has the permutation group $\Sym_3$ as a group of automorphisms. The
vertices of the diagram are labeled by the simple roots of the Lie
algebra of type $D_4$. Let $(e_1,e_2,e_3,e_4)$ be the standard
orthonormal basis of the Euclidean space $\mathbb R^4$. The simple
roots are
\[
\alpha_1= e_1 -e_2,\quad \alpha_2= e_2 -e_3,\quad \alpha_3= e_3 -e_4 \ \text{
  and} \ \alpha_4= e_3 +e_4
\]
(see \cite{bourbaki}).  The permutation $\alpha_1 \mapsto \alpha_4,\ 
\alpha_4 \mapsto \alpha_3, \ \alpha_3 \mapsto \alpha_1$, $\alpha_2
\mapsto \alpha_2 $ is an automorphism of order~$3$ of the Dynkin
diagram.  Its extension to a linear automorphism of $\mathbb R ^4 $ is
given by the orthogonal matrix $\mu$ of \eqref{matrixmu.equ}.
The matrix
\begin{equation}\label{matrixnu.equ}
\nu  =\left( \begin{matrix}
    1 &  0 &  0 & 0 \\
    0 & 1 &  0 &  0\\
    0 & 0 & 1 & 0 \\
    0 & 0 & 0 & -1
\end{matrix} \right)
\end{equation}
extends the automorphism $\alpha_1 \mapsto \alpha_1$, $\alpha_4 \mapsto
\alpha_3$, $\alpha_3 \mapsto \alpha_4$, $\alpha_2 \mapsto \alpha_2 $.
The set  $\{\mu, \nu\}$ generates a subgroup of $\OO_4$
isomorphic to $\Sym_3$, which restricts to the automorphism group of
the Dynkin diagram.  The group
\[ 
W(D_4) =(\Sym_2\wr \Sym_4)^+ =\Sym_2 ^3 \rtimes \Sym_4
\]
is the Weyl group of the split adjoint algebraic group $\PGO^+_8$,
which is of type $D_4$.


The group $W(D_4) =\Sym_2^{3} \rtimes \Sym_4 $, as a subgroup of the
orthogonal group~$\OO_4$, is generated by the reflections with respect
to the roots of the Lie algebra of $\PGO^+_{8}$.  Elements of $
\Sym_2\wr\Sym_4$ can be written as matrix products
  \begin{equation} \label{rem:matrixrepcoverings}
  w = D \cdot P(\pi),
  \end{equation}
  where $D$ is the diagonal matrix
  $\diag(\varepsilon_1,\varepsilon_2,\varepsilon_3, \varepsilon_4)$,  $\varepsilon_i = \pm 1$,
  and $P(\pi)$
  is the permutation matrix of $\pi \in \Sym_4$.  The group $
  \Sym_2\wr\Sym_4$ fits into the exact sequence
\begin{equation}\label{beta.equ}
1 \to \Sym_2^4 \to \Sym_2\wr\Sym_4 \overset{\beta}{\to} \Sym_4 \to 1
\end{equation}
where $\beta$ maps $w = D \cdot P(\pi)$ to $\pi$. Elements of $
W(D_4)$ have a similar representation, with the supplementary
condition $\prod_i \varepsilon_i =1$.

In relation with the geometric description of $C_1$ and $C_2$ at the
end of \S\ref{sec:trialcov}, note that the
group $\Sym_2\wr \Sym_4 = \Sym_2^4 \rtimes \Sym_4$ is the group
of automorphisms of the hypercube~$\mathcal K$. The subgroup $W(D_4)
=(\Sym_2\wr \Sym_4)^{+}= \Sym_2^3 \rtimes \Sym_4$ consists of the automorphisms
of $\mathcal K$ respecting the decomposition of the set of vertices
as $X \sqcup Y$, i.e., automorphisms of the half-hypercube.

\subsection{Automorphisms of $W(D_4)$}
We view $W(D_4)$ as a subgroup of $\OO_4$ as in
\eqref{rem:matrixrepcoverings}. Conjugation $x \mapsto \mu x \mu^{-1}$ with the matrices $\mu$
and $\nu$ on $\OO_4$ induce by restriction outer automorphisms $\tilde
\mu$ and $\tilde \nu$ of $W(D_4)$. The set $\{\tilde
\mu,\tilde \nu\}$ generates a group of automorphisms of $W(D_4)$
isomorphic to $\Sym_3$ (see already \cite[p.~368]{cartan}).
The center of $W(D_4)$ is isomorphic to $\Sym_2$ and is generated by
\begin{equation} \label{equ:center}
w_0 =\diag(-1,-1,-1,-1) = -1.
\end{equation}
Thus $ W(D_4)/\langle w_0\rangle $ acts on $ W(D_4)$ as the group of inner 
automorphisms. Let $\psi$ be the automorphism of order $2$ of $W(D_4)$
given by
 \[
 \psi \colon D \cdot P(\pi )\mapsto D \cdot P(\pi)\cdot (w_0
 )^{\sgn(\pi)},
\] 
or equivalently by $x \mapsto x\det(x)$, $ x \in
 W(D_4) \subset\OO_4 $.
 A proof of the following result can be found in
 \cite[Theorem~31,(5)]{franzsen_howlett} or in \cite[Prop.~2.8,(e)]{franzsen}:
\begin{prop} \label{autWplus.prop}
 \[
 \Aut\big(W(D_4)\big) \simeq \big( (W(D_4)/\langle w_0\rangle) \rtimes
 \,\Sym_3\big) \times \langle\psi\rangle.
\]
\end{prop}
For any $w\in W(D_4)$, we let $\Int(w)\colon x \mapsto wxw^{-1}$ be the inner automorphism of
$W(D_4)$ defined by conjugation by $w$, and by
$\Int\bigl(W(D_4)\bigr)$ the group of inner automorphisms of
$W(D_4)$. As an immediate consequence of
Proposition~\ref{autWplus.prop}, we have

\begin{cor}
\label{autWplus.cor}
\[
\Aut\bigl(W(D_4)\bigr)/\Int\bigl(W(D_4)\bigr) \simeq \Sym_3\times
\langle\psi\rangle.
\]
\end{cor}

We call \emph{trialitarian} the outer automorphisms of order~$3$ of
$W(D_4)$. As observed above, the automorphisms $\tilde\mu$ and $\tilde
\mu^2$ are trialitarian. Conjugation by the matrix $\rho$ of
\eqref{matrixrho.equ} also yields a trialitarian automorphism
$\tilde\rho$: indeed, we have $\rho^3=1$ and
\[
\rho \mu^{-1} = -1 \cdot\left( \begin{matrix}
    0 & 0 &  0&  1\\
    0&  1&   0 &  0\\
    1 &  0&   0 &  0\\
    0&   0&  1 &  0\\
\end{matrix}\right) \in W(D_4);
\]
hence, letting $w$ be the matrix on the right side, we
have $\tilde\rho=\Int(w)\circ\tilde\mu$.

\begin{prop}\label{prop:conjclasses}
  Any trialitarian automorphism of $W(D_4)$ is conjugate in the group
  $\Aut\big(W(D_4)\big)$ to either $\tilde \mu$ or $\tilde \rho$.
\end{prop}

\begin{proof}
  Explicit computation (with the help of the algebra computational
  system Magma \cite{magma}) shows that
  the conjugation class of $\tilde \rho$ contains $16$ elements and
  the conjugation class of $\tilde \mu$ contains $32$ elements. In
  view of Proposition~\ref{autWplus.prop} any trialitarian
  automorphism of $W(D_4)$ is the restriction to $W(D_4)$ of
  conjugation by an element $u\in\OO_4$ of the form $u=\mu\cdot w$ or
  $u=\mu^2\cdot w$, with $w\in W(D_4)$ and $u^3=1$. There are $48$
  elements $u$ of this form, hence the claim.
\end{proof}
\begin{cor}
There are up to isomorphism two types of subgroups of fixed points
of  trialitarian automorphisms of $W(D_4)$, those isomorphic to 
$\Fix(\tilde\mu)$  and those isomorphic to $\Fix(\tilde\rho)$ .
\end{cor}
\begin{proof}
Trialitarian automorphisms which are conjugate  in  $\Aut\big(W(D_4)\big)$ have isomorphic groups of fixed
points.
\end{proof}
\begin{prop} 
  \label{prop:fixmu}
  1)  The $2$-dimensional subspace of $\mathbb R^4$ generated by the set
  of elements $\{e_1 +e_3, e_2 - e_3\}$
  is fixed under $\mu$. \\
 2)  The set $\{e_1 +e_3, e_2 - e_3\}$ generates a root system of type
  $G_2$ and the group $\Fix(\tilde \mu)$ is the
  corresponding Weyl group, which is the dihedral group $\Dih_6$ of
  order $12$.
\end{prop}
 
\begin{proof}
  By explicit computation.
\end{proof}
\begin{prop}
  \label{prop:fixrho}
  The group $\Fix(\tilde \rho)$ is isomorphic to the group of order
  $24$ of Hurwitz quaternions $\mathbb H^1 = \Quat_8 \rtimes \Cyc_3$.
  This group is isomorphic to the double covering $\tilde \Alt_4$ of
  \,$\Alt_4$.
\end{prop}
\begin{proof}
  Recall that the matrix $\rho$ is obtained by choosing $(1,\, i,\,
  j,\, k)$ as 
  basis of $\mathbb R^4$ and letting $ -\frac{1}{2}(1 +i+j+k)$ operate
  by left multiplication in $\mathbb H$. The group $\mathbb H ^1$ has
  a representation in $W(D_4)$ by right multiplication which obviously
  commutes with the action of $\rho$. Hence $\Fix(\tilde \rho)$
  contains a copy of $\mathbb H ^1$. The claim then follows from the
  fact that $\Fix(\tilde \rho)$ has $24$ elements.
\end{proof}

\subsection{Cohomology with $W(D_4)$ coefficients}

Each automorphism $\alpha\in\Aut\bigl(W(D_4)\bigr)$ acts on
$H^1\bigl(\Gamma,W(D_4)\bigr)$ by
\[
\alpha_*\colon [\varphi] \mapsto [\alpha\circ\varphi],
\]
where $\varphi \colon \Gamma \to W(D_4)$ is a cocycle with values in
$W(D_4)$.
If $\alpha'=\Int(w)\circ\alpha$ for some $w\in W(D_4)$, then for all
cocycles $\varphi\colon\Gamma\to W(D_4)$ we have
\[
w\cdot\alpha\bigl(\varphi(\gamma)\bigr)\cdot w^{-1}
=\alpha'\bigl(\varphi(\gamma)\bigr) \qquad\text{for all
  $\gamma\in\Gamma$},
\]
hence $[\alpha\circ\varphi] = [\alpha'\circ\varphi]$ and therefore
$\alpha_*=\alpha'_*$. Thus, the action of $\Aut\bigl(W(D_4)\bigr)$ on
$H^1\bigl(\Gamma,W(D_4)\bigr)$ factors through
$\Aut\bigl(W(D_4)\bigr)/\Int\bigl(W(D_4)\bigr)\simeq \Sym_3\times
\langle\psi\rangle$. In particular the symmetric group $\Sym_3$
acts on $H^1\bigl(\Gamma,W(D_4)\bigr)$.
Under the bijections~\eqref{equ:etexcovplus}, the symmetric group $\Sym_3$
also acts on $\Iso\big((\Cov_\Gamma^{2/4})^+\big)$ and 
$ \Iso\bigl((\Etex^{2/4}_\Gamma)^+\bigr)$. The action of the
outer automorphism $\tilde \nu$ associates to the oriented 
$2/4$-covering $(Z/Z_0, \partial_Z)$ the oriented covering
$\kappa(Z/Z_0,\partial) = (Z/Z_0,  \partial_Z \circ \iota) $ where
$\boldsymbol{2}\xleftarrow{\iota}
\boldsymbol{2}$ twists the orientation. The proof of
Theorem~\ref{thm:Candtriality} shows that the action of $\tilde\mu$
maps the class of an oriented covering $(Z/Z_0,\partial_Z)$ to the
class of $C_1^+(Z/Z_0,\partial_Z)$.  

\section{Triality and \'etale algebras}
We next
investigate triality  on isomorphism classes of \'etale algebras 
using Galois cohomology.
Oriented extensions of \'etale algebras $S/S_0$ with $\dim_FS=8$ and 
$\dim_FS_0=4$  correspond to cocycles, i.e., continuous homomorphisms $\Gamma \to W(D_4)$, and isomorphism classes of such algebras correspond to
cocycles up to conjugation. 
If the cocycle factors through a subgroup $G$ of $W(D_4)$,  the 
conjugacy class of $G$ in $W(D_4)$ is determined by the isomorphism class
of the algebra. Thus it makes sense to classify isomorphism classes of
algebras according to the conjugacy classes
of the subgroups $G$  of $W(D_4)$. 
 
We give in Table~1 a list of all  conjugacy classes of 
subgroups of  $W(D_4)$. We still consider  $W(D_4)$ as a subgroup of $\OO_4$ 
(see \eqref{rem:matrixrepcoverings}) and use the following notation.
The group $W(D_4)$ fits into the split exact sequence:
\begin{equation} \label{equ:Wplussequence}
1 \to \Sym_2^3 \to W(D_4) \overset{\beta}{\to} \Sym_4 \to 1
\end{equation}
where $\beta$ is as in \eqref{beta.equ}.  For each
subgroup $G$ of $W(D_4)$ we denote by $G_1$ the restriction $G \cap
\Sym_2^3 $ and by $G_0$ the projection $\beta(G)$.  The center of
$W(D_4)$, generated by $w_0 =\diag(-1,-1,-1,-1) = -1$ is denoted by
$C$ and we set $ w_1=\diag(1,-1,1,-1),\,w_2=\diag(1,-1,-1,1)$ and $w_3
=\diag(-1,-1,1,1$) for special elements of the subgroup $\Sym_2^3
\subset W(D_4)$ given by diagonal matrices. We denote by $\Sym_n$ the permutation group of $n$
elements, $\Alt_n$ is the alternating subgroup, $\Cyc_n$ is cyclic of
order $n$, $\Dih_n$ is the dihedral group of order $2n $, $\Vier_4$ is
the Klein $4$-group, and $\Quat_8$ is the quaternionic group
with eight elements. We refer to \cite{CHMcK} for a description of the groups $[2^2]4$ and $\Quat_8:2$ in Table 1.  In Column $S$ we summarize the various
possibilities for \'etale algebras of dimension $8$ associated to the
class of a cocycle $\alpha\colon \Gamma \to W$ which factors through $G$
and in Column $S_0$
\'etale algebras of dimension $4$ associated to the class of the
induced cocycle $\beta\circ \alpha \colon \Gamma \to \Sym_4$
which factors through $G_0$. The entry
$K$ in one of the columns $S$ or $S_0$ denotes a quadratic separable field
extension.
We use
symbols $E$, respectively $E_0$ for separable field extensions  whose
 Galois closures have Galois groups $G$,
respectively $G_0$.  The symbol $R(E)$ stands for the cubic
resolvent of $E$ if $E$ is a quartic separable field extension and 
$\lambda^2E$ stands for the second lambda power of $E$
(see \cite{knustignol}, where it is denoted $\Lambda_2(E)$ or
\cite{GMS}, where it is denoted $E(2)$\footnote{The corresponding
  $\Gamma$-set $\lambda^2X$ is obtained by removing the diagonal from
  $X \times X$ and dividing by the involution $(x,y) \mapsto
  (y,x)$.}). If $\dim_FE = 4$,  $\lambda^2E$ admits an involution $\sigma$
  and for any quadratic \'etale algebra $K$ with involution $\iota$ we set 
  $K * \lambda^2E= (K \tens_F \lambda^2E) ^{\iota \tens \sigma}$. 
  We write $\overline{E}_0$ for the
Galois closure of $E_0$.  The symbol $\ell$ gives the
number of subgroups in the conjugacy class of $G$, MS refers to the
maximal subgroups of $G$ and in column $T$ we give the two conjugacy
classes which are the images of the class of $G$ under the
trialitarian automorphisms $\tilde \mu$ and $\tilde \mu^2$.

Entries N, $|G|,\, \ell$ and MS in the table
were generated with the help of the Magma algebra 
software~\cite{magma}. The computation of the entry~$T$, the explicit
representation of the group $G$ as an exact sequence and the decomposition of
the \'etale algebras as products of fields were checked case by case.

Explicit computations of trialitarian triples were made in
\cite{beltrametti} and \cite{zweifel} using the description
of the trialitarian action given in the proof of
Theorem~\ref{thm:Candtriality}.

\vfill

\vfill
 
 \begin{center}

\begin{minipage}[t]{15cm}
  \footnotesize


\[
\begin{array}{|l|l|l|l|l|l|l|l|l|l|l|}
\hline
N & S_0&S &  G_1 \to G\to  G_0^{\phantom{S^S}} & \vert G\vert & \ell& MS & T\\ 
\hline
1  & F^4 & F^8 &1 \to 1   \to  1 ^{\phantom{S^S}}                  & 1 & 1& 1   & 1, 1 \\
2  &  F^4 & K^4 & C \to \Sym_2 \to 1                    & 2& 6 &1 &  2,2 \\
3 & K^2& K^4&  1 \to \Sym_2 \to \Sym_2 \subset \Vier_4      &2 &6 &  1&  3,5 \\
4 &  K^2& K^4&  1 \to \Sym_2 \to \Sym_2\subset \Vier_4       &2 &6 &  1&  4,3 \\
5 & F^4 & F^2 \times K^2 & \langle w_1\rangle \to \Sym_2 \to 1             &2 &6 &  1&  5,4 \\
6& F^2 \times K  & F^2 \times K^2 & 1 \to \Sym_2 \to \Sym_2 \not \subset \Vier_4  &2 &12 &  1&  6,6 \\
7 &  F^2 \times K  & F^2 \times K^2 & 1 \to \Sym_2 \to \Sym_2 \not \subset \Vier_4  &2 &12 &  1&  7,7 \\
8 & F \times E_0&  F^2 \times S^2_0 & 1\to \Cyc_3 \to \Cyc_3                   &3 &16&  1&  8,8\\
9 & F^4 & K_1 ^2 \times K_2^2 & \langle C,w_1\rangle \to \Sym_2^2 \to 1               & 4&3&  2\ 5 & 11,10\\
10 & K ^2 & K \tens K_1 ^2  & C \to \Sym_2^2 \to  \Sym_2 \subset \Vier_4            &4&3  & 2\ 3 & 9, 11\\

11   & K ^2 & K \tens K_1 ^2  & C \to \Sym_2^2 \to  \Sym_2 \subset \Vier_4             &4 & 3 & 2\ 4 & 10,9\\
12 & K_1 \tens K_2 &  K_1 \tens K_2^2  & 1 \to \Sym_2^2 \to   \Vier_4              & 4& 4 &  4&14,13\\
13 & F^4 & K_1^2 \times K_2^2 & \langle w_1,w_2\rangle \to \Sym_2^2 \to 1                & 4 &4 & 5 & 12,14 \\
14  & K_1 \tens K_2 &  K_1 \tens K_2^2   &  1 \to \Sym_2^2 \to   \Vier_4              &  4&4 & 3 & 13,12\\
15  & F^2 \times K& F^4 \times K_1 \tens K& \langle w_1\rangle \to \Sym_2^2 \to \Sym_2 \not \subset \Vier_4  &4 &6 & 5 \ 6&  18,17 \\
16  & K_1 \times K_2 & K_1 \tens K_2^2 & 1 \to \Sym_2^2 \to \Sym_2^2  &4 &6 & 4 \ 7&  21,19 \\
17  & K_1 \times K_2 & K_1^2 \times K_2^2 & 1 \to \Sym_2^2 \to \Sym_2^2   &4 &6 & 3\ 6&  15, 18\\
18  & K_1 \times K_2 & K_1^2 \times K_2^2 & 1 \to \Sym_2^2 \to \Sym_2^2   &4 &6 & 4 \ 6&  17, 15 \\
19  & F^2 \times K &K_1^2 \times K_1 \times K & \langle w_1\rangle \to \Sym_2^2 \to \Sym_2\not \subset  \Vier_4  &4 &6 & 5\ 7&  16,21 \\
20& K ^2& E^2 & C \to \Cyc_4 \to \Sym_2 \subset \Vier_4  & 4&2 &2& 20,20\\
21 & K_1 \times K_2 & K_1\tens K_2^2&  1 \to \Sym_2^2 \to \Sym_2^2   &4 &6 & 3\ 7& 19, 16 \\
22 & K_1 \times K_2 & K_1^ 2 \times K_1\tens K_2&  1 \to \Sym_2^2 \to \Sym_2^2  &4 &12 &  4\ 6\ 7& 23, 27 \\
23 &  K_1 \times K_2& K_1^2 \times K_1 \tens K_2&  1 \to \Sym_2^2 \to \Sym_2^2    &4 &12 &3 \ 6\ 7& 27, 22 \\
24 &  K ^2 & K^2 \times K \tens K_1&  \langle w_1\rangle \to \Sym_2^2 \to \Sym_2\subset  \Vier_4  &4 &12 & 3\ 4\ 5&  24,24 \\
25 & F^2 \times K & F^4 \times E &  \langle w_1\rangle \to \Cyc_4\to \Sym_2 \not \subset \Vier_4  &4 &12 &5 & 26, 28 \\
26 &E_0 & E_0^2& 1 \to \Cyc_4 \to \Cyc_4                    & 4 &12 & 4 & 28,25 \\
27 & F^2 \times K & K_1^2 \times K_1 \tens K & \langle -w_1\rangle \to \Sym_2^2 \to \Sym_2  \not \subset \Vier_4   & 4&12 &5\ 6\ 7 & 22,23 \\
28 & E_0 & E_0 \times  E_0  & 1  \to \Cyc_4 \to \Cyc_4                 & 4 &12 & 3 & 25, 26 \\
29 & F^2 \times K & K^2 \times K \tens K_1  &  C \to \Sym_2^2 \to \Sym_2\not \subset  \Vier_4  &4 &12 & 2\ 6\ 7& 29,29\\
30 &  F\times E_0 & F^2 \times E_0 \tens \Delta(E_0) &  1 \to \Sym_3 \to \Sym_3    & 6 &16 & 7\ 8& 30,30\\
31 & F \times E_0 & F^2 \times E &  C \to \Cyc_6 \to \Cyc_3         & 6 & 16 &   2\  8& 31,31 \\
32  &  F \times E_0 &  F^2 \times E_0^2 &  1 \to \Sym_3 \to \Sym_3      & 6 &16 & 6\ 8& 32,32\\
33  & E_0 & E & C \to \Sym_2^3 \to \Vier_4          & 8 & 1 & 11\ 12 & 34,35 \\
34 & E_0 & E  & C \to \Sym_2^3 \to \Vier_4                  & 8 & 1 & 10\ 14& 35,33 \\
35 & F^4 & K_1 \times K_2 \times K_3 \times K_{123} & \Sym_2^3 \to \Sym_2^3 \to 1                & 8 & 1 & 9\ 13& 33, 34 \\
36 & E_0 & E &  C \to \Quat_8 \to \Vier_4                   & 8& 2 & 20 & 36,36\\
37 & K^2 & E_1 \times E_2 & \langle C,w_1\rangle \to \Sym_2 \times \Cyc_4 \to \Sym_2 \subset  \Vier_4  &8&3 & 3\ 20 & 38,40 \\
38 & E_0 & E_0 \tens K& C \to \Sym_2 \times \Cyc_4 \to \Cyc_4 &8&3 & 11 \ 20&  40,37 \\
39 & K^2 & K\tens K_1 \times K\tens K_2 &  \langle C,w_1\rangle \to \Sym_2 ^3 \to \Sym_2 \subset \Vier_4  &8& 3& 9\ 10\ 11 \ 24 & 39,39 \\
40 & E_0 & E_0 \tens K& C \to \Sym_2 \times \Cyc_4 \to \Cyc_4  &8&3  & 10 \ 20& 37,38\\
41 &K ^2K & E^2  & \langle C,w_2\rangle \to \Dih_4 \to \Sym_2 \subset  \Vier_4  & 8 & 6 & 9\ 11\ 20 & 49,45\\
42 & F^2 \times K& K_1 \times E  & \langle C,w_1\rangle \to \Sym_2 \times \Cyc_4 \to \Sym_2 \not \subset  \Vier_4  &8&6 & 9\ 25 & 44,47\\
43 & K_1 \times K_2 & K \tens K_1 \times K \tens K_2 & C \to \Sym_2^3 \to  \Sym_2^2   & 8&6& 10\ 17 \ 21\ 23\ 29 & 48,46\\
44 & E_0 & E & C \to \Sym_2 \times \Cyc_4 \to \Cyc_4& 8 & 6 & 11\ 26 & 47,42\\  
45 & K ^2 &  E ^2 & \langle C, w_2\rangle \to \Dih_4 \to \Sym_2 \subset  \Vier_4  & 8 & 6 & 9\ 10\ 20& 41,49\\
46& K_1 \times K_2 & K \tens K_1 \times K\tens K_2 & C \to \Sym_2^3 \to \Sym_2^2  &8&6 & 11\ 16 \ 18 \ 22\ 29&43,48\\
47& E_0 & E_0 \tens K&  C \to \Sym_2 \times \Cyc_4 \to \Cyc_4  &8&6& 10 \ 28& 42, 44\\
48 & F^2 \times K & K_1^2 \times K \tens K_2 & \langle C,w_1\rangle \to \Sym_2^3 \to \Sym_2 \not \subset  \Vier_4  &8&6& 9 \ 15\ 19\ 27\ 29 &46,43\\
49 & E_0 & E & C \to \Dih_4 \to \Vier_4 & 8&6& 10\ 11\ 20& 45,41\\
\hline
\end{array}
\]
\end{minipage}\\[2ex]

Table 1
\end{center}
\begin{center}

\begin{minipage}[t]{15cm}
  \footnotesize
\[
\begin{array}{|l|l|l|l|l|l|l|l|l|l||l|}
\hline
N & S_0 & S& G_1 \to  G \to   G_0^{\phantom{S^S}} & \vert G\vert & \ell & MS & T\\ 

\hline

50 & F^2 \times K & K^2 \times E^{\phantom{S^S}} & \langle w_1,w_2\rangle \to \Dih_4 \to \Sym_2 \not \subset  \Vier_4  &8&6& 13\ 19\ 25&55,51\\
51 & E_0 & E & 1 \to \Dih_4 \to \Dih_4 & 8 & 12 & 14\ 21\ 2 8 &50,55\\
52 & E_0 & E_0^2 & 1 \to \Dih_4 \to \Dih_4 & 8 & 12 & 14\ 17\ 28 & 54, 57 \\
53 & K_1 \times K_2 & K\tens K_1 \times K\tens K_2& \langle w_1\rangle \to \Sym_2^3 \to   \Sym_2^2  &8&12&16\ 19\ 21\ 22\ 23\ 24\ 27 &53,53\\
54 & F^2 \times K & F^2 \times K \times E  & \langle w_1,w_2\rangle \to \Dih_4 \to  \Sym_2 \not \subset  \Vier_4  & 8& 12 &13\ 15\ 25&57,52\\
55 & E_0 & \bar E_0 & 1 \to \Dih_4 \to \Dih_4 &  8& 12 &12\ 16\ 26&51,50\\
56 & K_1 \times K_2 & K_1^2 \times E & \langle w_1\rangle \to \Sym_2^3 \to \Sym_2^2   &8&12 & 15\ 17\ 18\ 22\ 23\ 24\ 27& 56,56\\
57&  E_0 & E_0^2 & 1 \to \Dih_4 \to \Dih_4 &  8& 12 &12\ 18\ 26 & 52,54\\
58 & E_0 & E_0^2 & 1 \to \Alt_4 \to \Alt_4 &12 & 4 & 8\ 14 & 59,60\\
59 & F \times E_0 & F^2 \times E_0^2 & \langle w_1,w_2\rangle \to \Sym_2^2 \rtimes \Cyc_3 \to \Cyc_3 &12 & 4 & 8 \ 13& 60,58\\
60 & E_0 & E_0 ^2 & 1 \to \Alt_4 \to \Alt_4 &12 &4 & 8\ 12 & 58,59 \\
61 & F \times E_0 & K \times K \tens E_0& C \to \Sym_2 \times \Sym_3 \to \Sym_3 &12 & 4 &29, 30 ,31,32& 61,61\\
62 & K_1 \times K_2 & E^2 & \langle C,w_1\rangle \to[2^2]4 \to \Sym_2^2   &16 & 3 & 39\ 42 & 66, 63 \\
63 & E_0 & E & \langle C,w_1\rangle \to [2^2]4 \to  \Cyc_4  &16 & 3 & 39\ 47 & 62, 66 \\
64 & K ^2& E_1 \times E_2& \Sym_2^3 \to \Sym_2 \times \Dih_4 \to \Sym_2 \subset  \Vier_4  & 16& 3 & 35\ 37\ 39\ 41\ 45& 68,67\\
65 & K_1 \times K_2 & K_1 \tens K_3 \times K_2 \tens K_4 & \langle C,w_1\rangle \to \Sym_2^4  \to \Sym_2^2 & 16 & 3 & 39\ 43\ 46\ 48\ 53\ 56& 65,65\\
66 & E_0 & E & \langle C,w_1\rangle \to [2^2]4 \to \Cyc_4 & 16 & 3 & 39\ 44 & 63,62 \\
67 & E_0 & E & \langle C,w_1\rangle \to \Sym_2 \times \Dih_4 \to \Vier_4 & 16 & 3 & 34\ 39\ 40\ 45\ 49& 64,68\\
68 & E_0 &E & \langle C,w_2\rangle \to \Sym_2 \times \Dih_4 \to \Vier_4 & 16 & 3 & 33\ 38\ 39\ 41\ 49 & 67,64\\
69& K_1 \times K_2 & E^2 &\langle C,w_1\rangle \to [2^2]4 \to \Sym_2^2  &16 & 6 & 37\ 42\ 48& 73,72 \\
70 & E_0 & E & \langle C,w_1\rangle \to \Quat_8:2 \to \Vier_4 & 16 &6 &36\ 37\ 38\ 40\ 41\ 45\ 49& 70,70\\
71 & F^2 \times K & K_1^2 \times E & \Sym_2^3 \to \Sym_2 \times \Dih_4 \to \Sym_2 \not \subset  \Vier_4  &16 & 6 & 35\ 42\ 48\ 50\ 54&75,74\\
72 & E_0 & E & \langle C,w_1\rangle  \to [2^2]4 \to \Vier_4 &16 & 6& 40\ 43\ 47& 69,73 \\ 
73 & E_0 & E & C \ \to[2^2]4 \to \Dih_4 &16 & 6& 38\ 44\ 46& 72,69\\ 
74 & E_0 & E_0 \tens K & C \to \Sym_2 \times \Dih_4 \to \Dih_4 & 16 &6 &34\ 43\ 47\ 51\ 52& 71,75\\
75 & E_0 & E_0 \tens K & C \to  \Sym_2 \times \Dih_4 \to \Dih_4 & 16 & 6& 33\ 44\ 46\ 55\ 57&74,71\\
76 & E_0 & E_0^2 & 1 \to \Sym_4 \to \Sym_4 &24 & 4 & 32\ 52\ 58& 79,81\\
77& F \times R(E) & \Delta(E) \times \lambda^2E& \langle w_1,w_3\rangle \to \Sym_2^2 \rtimes \Sym_3 \to \Sym_3 &24&4& 30\ 50 \ 59& 84,82\\
78 & F \times E_0 & \Delta(E) \times \lambda^2E & \Sym_2^3 \to \Sym_2^3 \rtimes \Cyc_3 \to \Cyc_3 & 24 & 4 & 31\ 35\ 59\  & 83,80 \\
79 & F \times R(E) & F^2 \times \lambda^2E &  \langle w_1,w_3\rangle \to \Sym_2^2 \rtimes \Sym_3 \to \Sym_3 &24&4& 32\ 54\ 59 & 76,81\\
80 & E_0 & E_0 \tens K & C \to \Sym_2 \times \Alt_4 \to \Alt_4 & 24 & 4 &31\ 34\ 58 & 78,83\\
81 & E_0 & E_0^2 & 1 \to \Sym_4 \to \Sym_4 &24 &4 &32\ 57\ 60& 76,79\\
82 & E_0 &E_0 \tens \Delta(E_0)  &  1 \to \Sym_4 \to \Sym_4 &24 &4 &30\ 51\ 58& 77,84\\
83 & E_0 & E & C \to \Sym_2 \times \Alt_4 \to \Alt_4 & 24 & 4 &31\ 33\ 60 & 78, 80\\
84 & E_0 & E_0\tens \Delta(E_0) &  1 \to \Sym_4 \to \Sym_4 &24 &4 &30\ 55\ 66& 82,77\\
85 & E_0  & E & \Sym_2 \to \tilde \Alt_4 \to \Alt_4 &24 & 4 &31 \ 36 & 85,85\\
86 & E_0 & E& \Sym_2^3 \to \Sym_2^3 \rtimes \Vier_4 \to \Vier_4 &32 & 1 & 64\ 67\ 68\ 70 &86,86\\
87 & E_0  & E  & \langle C,w_1\rangle \to \Sym_2^2 \rtimes \Dih_4 \to \Dih_4&32 & 3& 63\ 65\ 67\ 72\ 74 & 92,90\\
88 & E_0  & E &  \Sym_2^3 \to \Sym_2^3 \rtimes \Cyc_4 \to \Cyc_4 &32 & 3& 63\ 64 \ 66& 91,89\\
89&E_0  & E   & \langle C,w_1\rangle \to \Sym_2^3 \rtimes \Cyc_4 \to  \Dih_4&32 & 3&62\ 66\ 67& 88,91\\
90 & E_0  & E &  \langle C,w_1\rangle \to \Sym_2^2 \rtimes \Dih_4  \to \Dih_4 & 32 & 3 &65\ 66\ 68\ 73\ 75 & 92,87\\
91 & E_0  & E& \langle C,w_1\rangle \to \Sym_2^3 \rtimes \Cyc_4 \to  \Dih_4&32 & 3& 62\ 63\ 68& 89,88\\
92 & K_1 \times K_2 & E_1 \times E_2& \Sym_2^3 \to \Sym_2^2 \rtimes \Dih_4 \to \Sym_2^2   &32 &3 & 62\ 64\ 65\ 69\ 71&87,90\\
93 &  F \times R(E) & \Delta(E) \times K *\lambda^2E  & \Sym_2^3 \to  \Sym_2  \times \Sym_4 \to \Sym_3 & 48 &4 & 61\ 71\ 77\ 78\ 79& 94,95\\
94 & E_0 & E_0 \tens K & C \to  \Sym_2  \times \Sym_4 \to \Sym_4 & 48 & 4& 61\ 75\ 81\ 83\ 84& 95,93 \\
95 &E_0 & E_0\tens K & C \to  \Sym_2  \times \Sym_4 \to \Sym_4 & 48& 4& 61\ 74\ 76\ 80\ 82& 93,94\\ 
96 & E_0 & E & \Sym_2^3 \to \Sym_2^3 \rtimes \Dih_4 \to \Dih_4 &64 & 3& 86\ 87\ 88 \ 89 \ 90\ 91\ 92 & 96, 96\\
97 & E_0 & E & \Sym_2^3 \to \Sym_2^3\rtimes \Alt_4 \to \Alt_4 &96 &1& 78\ 80\ 83\ 85\ 86& 97, 97\\
98 &E_0 & E  &  \Sym_2^3 \to \Sym_2^3\rtimes \Sym_4 \to \Sym_4 &192 & 1 & 93\ 94\ 95\ 96\ 97 &98,98\\
\hline
\end{array}
\]
\end{minipage}\\[2ex]
Table 1
 \end{center}

 \subsection{Trialitarian triples and fixed points}
 Let $\alpha$ be any trialitarian automorphism of $W(D_4)$.
The subset $H^1\bigl(\Gamma,W(D_4)\bigr)^{\Cyc_3}$
of cohomology classes that are fixed under $\alpha_*$ is independent 
of the particular choice of $\alpha$. Clearly, the image in
$H^1\bigl(\Gamma,W(D_4)\bigr)$ of any cohomology class in
$H^1\bigl(\Gamma,\Fix(\alpha)\bigr)$ is fixed under $\alpha_*$, hence
this image lies in $H^1\bigl(\Gamma,W(D_4)\bigr)^{\Cyc_3}$. Thus, we have
a canonical map

\begin{equation} \label{eq:fixedcohomology}
H^1\bigl(F, \Fix(\alpha)\bigr) \to H^1\big(F,W(D_4)\big)^{\Cyc_3}
\end{equation}
for any trialitarian automorphism $\alpha$ of $W(D_4)$.

\begin{thm}
Any class in $H^1\big(F,W(D_4)\big)^{\Cyc_3}$ lies in the image of the
  map \eqref{eq:fixedcohomology} for $\alpha =\tilde \mu$ or
  $\alpha=\tilde \rho$ as in \eqref{prop:conjclasses}.
\end{thm}

\begin{proof} 
If the class $[\varphi]$ of a cocycle $\varphi\colon \Gamma \to W(D_4)$
belongs to  $H^1\bigl(\Gamma,W(D_4)\bigr)^{\Cyc_3}$, the cocycle factors through
a subgroup $G$ whose conjugacy class is invariant under~$\tilde \mu$.
 We get from Columns N and T of Table 1 a list of all triples of conjugate classes of subgroups  $W(D_4) $ which are permuted by a trialitarian 
 automorphism of $W(D_4)$ . 


A triple consists of three identic labels (for example $(70,70,70)$) if 
 there exists $a \in W(D_4)$ such that
 \[
 \mu G \mu^{-1} = a G a^{-1},
 \] 
  where $\mu$ is as in~\eqref{matrixmu.equ}.  A cocycle factoring through
  such a group $G$ does not necessarily correspond to  a triple  of  isomorphism classes of  \'etale algebras fixed under
 triality.  For a triple consisting of isomorphic algebras  we must have
 $a \in W(D_4)$ such that
 \begin{equation}\label{isotriples:equ}
 \mu x \mu^{-1} = axa^{-1}
 \end{equation}
 for all $x \in G$, since isomorphic classes of algebras are given by homomorphisms $\gamma\colon \Gamma \to G$ up to conjugation. Thus
 a necessary condition to get a triple of fixed isomorphism classes of algebras
 is that the conjugacy classes of all subgroups $H$ of $G$ are invariant
 under triality. Thus the conjugacy classes N = 24, 39, 53,  56, 65, 70, 86, 96, 97
 and 98 do not correspond to  triples of 
isomorphism classes of algebras invariant under triality. The conjugacy
classes left over in Table 1 are N= 1, 2, 6, 7, 8,
 20, 29, 30, 31, 32, 36, 61 and 85. They  give rise to triples of  isomorphism classes of \'etale 
 algebras fixed under
 triality, since they correspond to subgroups contained in $\Fix(\tilde \mu)$ 
(N = 61) or contained in  $\Fix(\tilde \rho)$ (N=85). 
\end{proof}

Observe that fixed \'etale algebras in class $N=61$ are of the form $K
\times (E_0 \tens K)$, where $K$ is quadratic and $E_0$ is cubic.
Hence they are not fields over $F$, in contrast to algebras in class
$N=85$.

\section{Trialitarian resolvents}

Trialitarian triples of \'etale algebras can be viewed as one \'etale
algebra with two attached resolvents (see Remark~\ref{resolvents:rem}). 
For example, let $E$ be a
quartic separable field with Galois group $\Sym_4$. The field $ E
\tens \Delta(E) $ is octic with the same Galois group $\Sym_4$ and the
extension $E\tens \Delta(E)/E$ corresponds to Class $N=82$ in Table 1.
Class $N=77$ in the same trialitarian triple corresponds to the
extension
\[
\big(\Delta(E) \times \lambda^2E\big)/\big(F \times R(E)\big),
\]
where $\Delta(E)$ is the discriminant, $R(E)$ is the cubic resolvent of $E$ and $\lambda^2E$ is the
second lambda power of $E$, as defined in \cite{knustignol}.

In this section we consider the situation where one
\'etale algebra in the triple is given by a separable polynomial and compute polynomials
for the two other  \'etale algebras. We  assume that the base field $F$ is infinite 
and has characteristic different from~$2$.
\begin{prop} \label{prop:primitiveelement}
  Let $S/S_0$ be an  \'etale algebra with involution of dimension $2n$ over $F$.   \\
  1) There exists an invertible element
  $x \in S$ such that $x$ generates $S$ and $x^2$ generates $S_0$. \\
  2) There exists a polynomial
\[
f_n(x) = x^n + a_{n-1}x^{n-1} + \dots + a_1x +a_0.
\]
with coefficients in $F$ such that $S_0 \iso F[x]/\big(f_n(x)\big)$
and
$S \iso F[x]/\big(f_{2n}(x)\big)$, where $f_{2n}(x) =f_n(x^2)$.\\
3) The algebra $S$ has trivial discriminant if and only if $(-1)^n a_0$ is a
square in $F$.
\end{prop}
\begin{proof}
  To prove 1) we are looking for
   invertible elements $x$ of $S$ such that $\Tr_{S/S_0}(x) = 0$ and such that the
  discriminant of the characteristic polynomial of $x^2$ is not zero. Any 
  such element generates $S$ and $x^2 \in S_0$ generates $S_0$. 
  These elements form an Zariski open subset of the space of trace
  zero elements.  One checks that this open subspace is not empty
  by going to an algebraic closure of $F$.\\
  2) follows from 1) and 3) follows from a discriminant formula (see
  \cite[p.~51]{brillhart}) for the discriminant $D(f_{2n})$ of $f_{2n}$:
\[
D(f_{2n}) = (-1)^n a_0 \cdot \big(2^nD(f_n)\big)^2
\]
(recall that $\Delta(S) \iso F[x]/\big(x^2 -D(f_{2n})\big)$).
 \end{proof} 
   
  \begin{thm}\label{prop:resolvents}
   Let $S/S_0$, $\dim_FS$= $8$,  with trivial discriminant,  be given as in
   Proposition~\ref{prop:primitiveelement}, by a polynomial
 \begin{equation}\label{eq:primitiveelement}
 f_4(x) = x^4 + ax^3 + bx^2 + cx +e^2.
 \end{equation}
     The polynomials
    \[ 
  \begin{array}{lll}
  f_4'(x) & = &x^4 +ax^3 +(\frac{3}{8} a^2 -\frac{1}{2}b + 3e)x^2\, + \\
  && (\frac{1}{16}a^3 -\frac{1}{4}ab  + c +\frac{1}{2} ae)x 
     + \big(\frac{1}{16}a^2- \frac{1}{4}b  -\frac{1}{2}e\big)^2 \quad \text{and}\\[2ex]
  f_4''(x) & = &  x^4 +ax^3 +(\frac{3}{8} a^2 - \frac{1}{2}b - 3e)x^2\, +\\
 &&  (\frac{1}{16}a^3 -\frac{1}{4}ab + c -\frac{1}{2} ae)x
   + \big(\frac{1}{16}a^2 -\frac{1}{4}b  +\frac{1}{2}e\big)^2 
  \end{array}
  \]
define extensions of \'etale algebras $S'/S'_0$,  $S''/S_0''$ 
such that the isomorphism classes of  $S/S_0$, $S'/S'_0$ and $S''/S_0''$
are in triality.
\end{thm}
\begin{proof}
Let $\{y_1,y_2,y_3,y_4\}$ be the set of zeroes of $f_4$ in a
 separable closure $F_s$ of $F$.  The set $\{\pm x_i =\pm \sqrt y_i,\ 
 i=1, \ldots 4\}$ is the set of zeroes of $f_8$.     
    Let $\xi$ be the column vector
    $[x_1,x_2,x_3,x_4]^T$.     If $\varphi: \Gamma \to W(D_4) \subset \OO_4$
    is the cocycle corresponding to $S/S_0 $, the group $\varphi(
    \Gamma )$ permutes the elements $\pm x_i$ through left matrix
    multiplication on $\xi$. The cocycle corresponding to $S'/S_0' $
    is given by
  \[
  \varphi': \gamma \mapsto \mu \varphi(\gamma)\mu^{-1}, \gamma
  \in\Gamma,
  \]
  where $\mu$ is as in~\ref{matrixmu.equ}.
  Thus $\varphi'( \Gamma )$ permutes the components of $\pm\xi'=\pm
  \mu\xi = \pm[x'_1,x'_2,x'_3,x'_4]^T$ and $\{\pm x'_i,
  i=1,\ldots,4\}$ is the set of zeroes of $f_8'$. It follows that
  \[
  f'_8(x) =f'_4(x^2)= \prod_i(x -x'_i)(x+x'_i) = \prod_i(x
  ^2-{x'_i}^2).
  \]
  The $x_i'$ are the components of $\xi' = \mu \xi$. Thus the
  coefficients of $ f'_8(x)$ can be expressed as functions of the
  $x_i$. Using that the symmetric functions in the $x_i$ can be
  expressed as functions of the coefficients of $f_8$ one gets (for
  example with Magma \cite{magma}) the expression given in
  Proposition~\ref{prop:resolvents} for $f_i'$. Similar computations
  with $\mu^2$ instead of $\mu$ lead to the formula for $f_4''$.
   \end{proof}
   
\begin{remark} Observe that we move from $f_8'$ to $f_8''$ by replacing
  $e$ by $-e$, as it should be.
\end{remark} 

\section{Triality and  Witt invariants of \'etale algebras}
The results of this section were communicated to us by J-P.  Serre,
\cite{Serre09}.  They are based on results of \cite{Serre07} and \cite{GMS}.
Similar results can be obtained for cohomological invariants of \'etale algebras instead of Witt invariants. 
Let $k$ be a fixed base field of characteristic not $2$ and $F/k$ be a 
field extension. Let $WGr(F)$ be the Witt-Grothendieck ring
and $ {W}(F)$ the Witt ring of $F$, viewed as functors of $F$.
 We recall that elements of
${ WGr}(F)$ are formal differences  $q -q'$ of isomorphism classes of nonsingular
quadratic forms over $F$ and that the sum and product are those induced by the
orthogonal sum and the tensor product of quadratic forms. The Witt ring
 ${ W}(F)$ is the quotient of ${ WGr}(F)$ by the ideal  consisting of
integral multiples of the $2$-dimensional diagonal form $\langle1,-1\rangle  $. 

Some of the following considerations hold for oriented quadratic extensions
$S/S_0$ of \'etale algebras of arbitrary dimension. To simplify notation 
we assume from now on that $\dim_F S=8$.

Let
$(\Etex^{2/4})^+$ be the functor which associates to $F$ the set
$\Big(\Etex^{2/4}_F\big)^+ $ of isomorphism classes of oriented quadratic extensions
$S/S_0$ of \'etale algebras over $F$ such that $\dim_FS= 8$. 
 A  \emph{Witt invariant} on $(\Etex^{2/4})^+ $, more precisely on $W(D_4)$,
 is a map 
 \[
 H^1\big(F, W(D_4)\big) \to
{W}(F)
\]
for each $F/k$, subject to compatibility and specialization conditions 
(see \cite{GMS}).
The set of Witt invariants \[
\Inv \!\big(W(D_4),{ W}\big) =  
\Inv\!\big((\Etex^{2/4})^+,{W}\big)
\] is  a module over ${W}(k)$. 
The aim of this section is to describe this module and how triality acts on it.  
A main tool is the following splitting principle, which is 
a special case of a variant of the splitting principle
for \'etale algebras (see \cite[Theorem 24.9]{GMS}) and which can be proved following the same lines.

\begin{thm}\label{thm:splitting} If $ a \in \Inv\big((\Etex^{2/4})^+,{W}\big)$ satisfies $a(S/S_0) =0$ for
every product of two biquadratic algebras
 \[
S = F\big(\sqrt{x}, \sqrt y\big) \times F\big(\sqrt{z}, \sqrt{t}\big), \quad 
S_0 = F\big(\sqrt{xy}\big) \times F\big(\sqrt{zt}\big). 
\]
over every extension $F$ of $k$, then $a= 0$.
\end{thm}

Let $G$ be an elementary abelian subgroup of $W(D_4)$ of type $(2,2,2,2)$. It 
belongs to the conjugacy class $N=65$ in Table 1. Theorem~\ref{thm:splitting}
 can be restated in the following form:
\begin{thm}\label{thm:splitting2} 
The restriction map \[
\Res \colon  \Inv\!\big(W(D_4),{ W}\big) \to \Inv\!\big(G,{ W}\big)
\]
is injective.
\end{thm}
\begin{proof}
$G$-torsors correspond to products of two biquadratic algebras.
\end{proof}

  A construction of Witt invariants is through trace forms.
Let $S/S_0 \in (\Etex^{2/4})^+$  
and let $\sigma$ be the involution of $S$. We may associate two nonsingular
quadratic trace forms to the extension
$S/S_0$:
\[
\begin{array}{lllll}
Q(x) &= & {Q}_S(x)& = &\frac{1}{8}\Tr_{S/F}(x^2) \\
Q'(x) &= & Q'_S(x)& =& \frac{1}{8}\Tr_{S/F}\big(x\sigma(x)\big), \ x \in S \\ 
\end{array}
\]
The decomposition
\[
S = \SSym(S,\sigma) \oplus \Skew(S,\sigma) 
\]
leads to orthogonal decompositions
\[
Q= Q^+ \perp Q^-, \quad    Q'= Q^+ \perp -Q^-,
\]
 hence the forms $Q^+$ and $Q^-$ define two Witt
invariants attached to $S/S_0$. The \'etale algebras associated to
$S/S_0$ by triality lead to corresponding invariants. We introduce
following notations:
$S/S_0 =S_1/S_{0,1}$ and  $S_i/S_{0,i},\ i=2,3$,  for  the
associated \'etale algebras. We denote the corresponding Witt invariants
by  $Q^+_i = Q^+ _{S_i}$ and $Q^-_i = Q^- _{S_i}$, $i=1,\,2$ and $3$.

Another construction of Witt invariants is through orthogonal representations.
Let~$\OO_n$ be the orthogonal group of the $n$-dimensional form $\langle1,\ldots,1\rangle  $.
Quadratic forms  over $F$ of  dimension $n$ are classified by the cohomology set
$H^1(F, \OO_n)$. Thus any group homomorphism $W(D_4) \to \OO_n$ gives rise to a Witt invariant. 
In particular we 
get a Witt invariant $q$ associated with
the orthogonal representation $W(D_4) \to \OO_4$ described in
\eqref{rem:matrixrepcoverings}.
Moreover the  group
 $W(D_4)$ has three normal subgroups $H_i$ of type $(2,2,2)$ (i.e., isomorphic
to $\Sym_2^3$), corresponding to the classes $N=33,\,34,\, 35$ of Table 1. 
Since the factor groups are isomorphic to $\Sym_4$, the canonical
representation $\Sym_4 \to \OO_4$ through permutation matrices leads to three
Witt invariants $q_1,\, q_2,\,q_3$.

\begin{prop}\label{prop:wittinv}
1) The Witt invariant $q$ is invariant under triality and coincides with $Q^-_i$,
$i=1,2,3$.\\
2) We have  $q_i = {Q_i^+}$, $i=1,\,2,\,3,$  and the three invariants 
$q_1,\,q_2,\,q_3$ are 
permuted by triality.
\end{prop}
\begin{proof}
The fact that $q$ is invariant under triality follows from the fact that
triality acts on $W(D_4)$ by an inner automorphism of $\OO_4$. Moreover the 
trialitarian action on  $W(D_4)$ permutes the normal subgroups $H_i$, hence 
the invariants $q_i$. For the other claims we may assume by the 
splitting principle that $S_1$ is a product of two biquadratic algebras
\begin{equation} \label{equ:biquadratic}
S_1 = F\big(\sqrt{x}, \sqrt{y}\big) \times F\big(\sqrt{z}, \sqrt{t}\big), \quad 
S_{0,1} = F\big(\sqrt{xy}\big) \times F\big(\sqrt{zt}\big).
\end{equation}
An explicit computation, using for example the description of triality
given in the proof of Theorem~\ref{thm:Candtriality} (see 
\cite{beltrametti} and \cite{zweifel}) shows that one can make the following
identifications
\[
S_2 =F\big(\sqrt{x}, \sqrt{z}\big) \times F\big(\sqrt{y}, \sqrt{t}\big), \quad
S_{0,2} = F\big(\sqrt{xz}\big) \times F\big(\sqrt{yt}\big)
\]
and
\[
S_3 =F\big(\sqrt{x}, \sqrt{t}\big) \times F\big(\sqrt{y}, \sqrt{z}\big), \quad
S_{0,3} = F\big(\sqrt{xt}\big) \times F\big(\sqrt{yz}\big).
\]
Observe that, with this identification, the $3$-cycle $(y,z,t)$ permutes
cyclically the algebras $S_i/S_{0,i}$. 
We get
\[
Q^-_i =\, \langle x,y,z,t\rangle   
\]
for $i=1,\,2,\,3$, and
\[
  Q^+_1 =\, \langle 1,1,xy,zt\rangle   ,\  Q^+_2 = \,\langle 1,1,xz,yt\rangle   ,  \ Q^+_3 =\, \langle 1,1,xt,yz\rangle  .
\]
The equalities $q=Q^-_i$ and $ q_i = {Q_i^+}$ follow from the fact 
that the corresponding cocycles are conjugate in $\OO_4$.
\end{proof} 

Further basic invariants are the constant invariant  $\langle 1\rangle  $ and  the discriminant
\[
\langle d\rangle  \, = \disc(q) =\disc(q_i),\  i=1,\,2,\,3,
\]
which corresponds to the $1$-dimensional representation 
$\det\colon W(D_4) \to \OO_1 = {\pm 1}$.
Since $Q^+_i(1) =1$,  the quadratic forms $q_i =Q_i^+$ represent
$1$ and  one can replace them by $3$-dimensional invariants 
$\ell_i = (1) ^\perp \subset q_i,\, i=1,2,3$. 
%

In the following result $\lambda^2q$  denotes the second exterior power of the quadratic form~$q$
(see \cite{GMS}). If $q = \langle \alpha_1,\ldots, \alpha_n\rangle  $ is diagonal, then $\lambda^2q$ is the
$n(n-1)/2$-dimensional form 
$\lambda^2q =\, \langle \alpha_1\alpha_2, \ldots, \alpha_{n-1}\alpha_n \rangle  $.
\begin{thm} \label{thm:basisW}
1)
The ${W}(k)$-module  $\Inv\!\big(W(D_4), W\big) =
\Inv\!\big((\Etex^{2/4})^+,{W}\big)$ is  free over ${W}(k)$ with basis
\begin{equation}\label{equ:basisW}
\big(\langle 1\rangle  , \,\langle  d\rangle  ,\, q, \langle d\rangle  \cdot\, q, \,\ell_1,\ell_2 , \, \ell_3\big).
\end{equation}
2) The elements $\langle 1\rangle  , \,\langle  d\rangle  \,, q,\,\langle d\rangle  \cdot q$ are fixed under triality and the
elements $\ell_1,\ell_2 \ \text{and} \ \ell_3$ are permuted.\\
3) The following nonlinear relations hold among elements of \eqref{equ:basisW}:
\[
\begin{array}{ccccc}
\langle d\rangle   &= & \disc(q)\, = \,\disc(q_i),\\
\lambda^2q &= &\ell_1 +\ell_2 + \ell_3\,  - \langle 1,1,1\rangle  \\
\langle 1, d\rangle   \cdot \,q  & = & q \cdot(\ell_i -\langle 1\rangle),\, i=1,\,2, \,3. \\
\end{array}
\]
\end{thm}
\begin{proof}
1) The proof follows the pattern of the proof of \cite[Theorem 29.2]{GMS}.
Let $G$ be an elementary  subgroup of  $W(D_4)$ of type $(2,2,2,2)$, i.e. isomorphic 
to $\Sym_2^4$. An arbitrary element of $H^1(F,G)$ is given by a
$4$-tuple $(\alpha_1,\alpha_2,\alpha_3,\alpha_4)
\in \big(F^\times/F{^{\times 2}}\big)^4$. For $I$ a
subset of $\boldmath{4} =\{1,2,3,4\}$, we write $\alpha_I$ for the product of the 
$\alpha_i$ for $i \in I$.
By \cite[Theorem 27.15]{GMS} the set $\Inv(G,W)$ is a free $W(k)$-module with 
basis $(\alpha_I)_{I\subset \boldmath{4}}$. It then follows from Theorem~\ref{thm:splitting2}
 that the family of elements given
in \eqref{equ:basisW} is linearly independent over $W(k)$.
Let $a$ be an element of $\Inv\!\big((\Etex^{2/4})^+,{W}\big)$ and let $S_\alpha$ 
be the algebra~\eqref{equ:biquadratic} for   
$\alpha_1=x$,  $\alpha_2=y$,  $\alpha_3=z$, $\alpha_4=t$. The map~$\alpha \mapsto a(S_\alpha)$ 
is a Witt invariant of $G$, hence by \cite[Theorem 27.15]{GMS} can be uniquely written 
as a linear combination
\begin{equation} \label{equ:lincombinv}
\sum c_I \cdot \langle \alpha_I\rangle   \ \text{with} \ c_I \in W(k).
\end{equation}
The  claim will follow if we show that the invariant $\alpha$ is in fact a linear combination of the elements
given in~\eqref{equ:basisW}. 
By \cite[Prop. 13.2]{GMS}, the image of the  restriction map 
$\Inv\!\big(W(D_4), W\big) \to \Inv(G, W)$ is contained in the $W(k)$-submodule of $\Inv(G, W)$ fixed
by the normalizer $N= \Sym_2^3 \rtimes \Dih_4$ of $G$ in $W(D_4)$ (conjugacy class
$N=96$ in Table~1). The group $N$ acts on the set of isomorphisms classes of
algebras $S_i$ by acting in the 
obvious way on the  symbols
$\pm \sqrt{x},\, \pm  \sqrt{y}, \, \pm \sqrt{z}, \, \pm  \sqrt{t}$.  
It follows that $N$ acts trivially on this set of isomorphisms classes. This shows that
only linear combinations of elements in the family 
\[
\begin{array}{ll}
\mathcal B  = & \{ \langle 1 \rangle,\, \langle x \rangle+
\langle y \rangle+  \langle z \rangle+ \langle t \rangle,\, \langle xy \rangle+
\langle zt \rangle, \,  \langle xz \rangle+ \langle yt \rangle, 
\langle xt \rangle+ \langle yz \rangle, \\
& \qquad 
 \langle  xyz\rangle+\langle xyt \rangle +
\langle xzt \rangle+
\langle yzt \rangle, \,\langle xyzt \rangle\}
\end{array}
\]
can occur in the sum~\eqref{equ:lincombinv}.
The  family $\mathcal B$ and the 
family given in~\ref{equ:basisW} are equivalent bases. This
implies the first claim of Theorem~\ref{thm:basisW}.
 Claim 2) follows from Proposition~\ref{prop:wittinv} and 3) is easy to check for
a product of biquadratic extensions.
\end{proof}

\end{document}